\documentclass[11pt,a4paper]{article}

\usepackage{amsmath}
\usepackage{amsfonts}
\usepackage[english]{babel}
\usepackage[latin1]{inputenc}
\usepackage{fontenc}
\usepackage{euscript}
\usepackage{vmargin}
\usepackage{amssymb}
\usepackage{latexsym}
\usepackage{enumerate}
\usepackage{graphicx}
\usepackage{hyperref}

\usepackage[titletoc,toc,title]{appendix}
\makeatletter
\renewcommand\paragraph{\@startsection{paragraph}{4}{\z@}%
            {-2.5ex\@plus -1ex \@minus -.25ex}%
            {1.25ex \@plus .25ex}%
            {\normalfont\normalsize\bfseries}}
\makeatother
\def\rain{\to +\infty}

\def\N{{\rm I\kern-.20em N}}
\def\R{{\rm I\kern-.20em R}}
\def\indi{{1\kern-.20em\rm I}}

\def\bkR{{\rm I\kern-.17em R}}
\def\bkN{{\rm I\kern-.20em N}}

\newtheorem{teo}{Theorem}[section]
\newtheorem{nota}{Remark}
\newtheorem{coro}{Corollary}[section]
\newtheorem{prop}{Proposition}[section]

\newtheorem{ex}{Example}[section]

\newcommand{\pg}{\hspace{0.6cm}}

\newcommand{\bdem} {\begin{proof}}
\newcommand {\edem}{\hfill $\square$ \end {proof}}
\def\X{{\bf X}}

\begin{document}
\title{Methods for estimating the upcrossings index: improvements and comparison}
\author{A.P.
Martins \footnote{
E-mail:amartins@ubi.pt}\\
Departamento de Matem\'atica \\
Universidade da Beira Interior\\ Portugal \and J.R. Sebasti\~ao \footnote{
E-mail:jrenato@ipcb.pt}\\ Escola Superior de Gest\~ao\\
Instituto Polit\'ecnico de Castelo Branco\\ Portugal}
\date{}
\maketitle

\noindent {\bf Abstract:}  The upcrossings index $0\leq \eta\leq 1,$ a measure of the degree of local dependence in the upcrossings of a high level by a stationary process, plays, together with the extremal index $\theta,$ an important role in extreme events modelling. For stationary processes, verifying a long range dependence condition,  upcrossings of high thresholds in different blocks can be assumed asymptotically independent and therefore blocks estimators for the upcrossings index can be easily constructed using disjoint blocks. In this paper we focus on the estimation of the upcrossings index via the blocks method and properties such as consistency and asymptotic normality are studied. We also enlarge the family of runs estimators of $\eta$ and provide an empirical way of checking local dependence conditions that control the clustering of upcrossings to improve the estimates obtained with the runs method.

We compare the performance of a range of different estimators for $\eta$ and illustrate the methods using simulated data and financial data.
\vspace{0.3cm}


\vspace{0.3cm}



\section{Introduction}\setcounter{equation}{0}

\pg Extreme Value Theory  aims to predict occurrence of rare events but with disastrous impact, making an adequate estimation of the parameters related with such events of primordial importance in Statistics of Extremes. The extremal index $\theta\in (0,1]$ and the upcrossings index $\eta\in (0,1]$ play an important role when modelling extreme events. The knowledge of these parameters, $\theta$ and $\eta,$ entails in particular the understanding of the way in which exceedances and upcrossings of high levels, respectively,  cluster in time. They provide different but complementary information concerning the grouping characteristics of rare events.

Suppose that we have $n$ observations from a strictly stationary process $\X=\{X_n\}_{n\geq 1},$
with marginal distribution function $F$ and finite or infinite right endpoint $x_F=\sup\{x:\ F(x)<1\}.$
$\X$ is said to have extremal index $\theta\in [0,1]$ if for each $\tau>0$ there exists a sequence of thresholds $\{u_n^{(\tau)}\}_{n\geq 1}$ such that
\begin{equation}\label{indextremal}
n(1-F(u_n^{(\tau)}))\xrightarrow [n\rain]{} \tau\qquad \textrm{and}\qquad F_n(u_n^{(\tau)})=P(M_n\leq u_n^{(\tau)})\xrightarrow [n\rain]{}e^{-\theta \tau},
\end{equation}
where $M_n=\max\{X_1,\ldots,X_n\}$ (Leadbetter (1983) \cite{lead1}). This parameter is a measure of clustering tendency of extremes, more precisely $\theta^{-1}$ is the limiting mean cluster size in the point process of exceedance times over a high threshold, under an appropriate long range dependence condition (Hsing {\it{et al.}} (1988) \cite{hsing}). For independent and identically distributed (i.i.d.) processes, extremes are isolated and hence $\theta=1,$ whereas for stationary processes, the stronger the dependence the larger the clusters of exceedances and therefore the smaller $\theta.$

Many results in extreme value theory may be naturally discussed in terms of point processes. When looking at the point process of upcrossing times of a threshold $u_n$, that is
\begin{equation}\label{processocruz}
\widetilde{N}_n(u_n)(B)=\sum_{i=1}^n\indi_{\{X_i\leq u_n<X_{i+1}\}}\delta_{\frac{i}{n}}(B),\quad B\subset (0,1],
\end{equation}
where $\indi_{A}$ denotes the indicator of event $A$ and $\delta _{a}$
the unit mass at $a,$ Ferreira (2006) \cite{fer1} showed that, under a long range dependence condition, if the sequence of point processes of upcrossings $\widetilde{N}_n(u_n)$ converges in distribution (as a point process on
(0,1]), then the limit is necessarily a compound Poisson process and the Poisson rate of the limiting point process is $\eta\nu,$ when the limiting mean number of upcrossings of $u_n$ is $\nu>0.$  That is, there is a clustering of upcrossings of high levels, where the underlying Poisson values represent cluster positions and the multiplicities are the cluster sizes. In such a situation there appears the  upcrossings index $\eta$ which has an important role in obtaining the sizes of the clusters of upcrossings.

The stationary process $\X$ is said to have upcrossings index $\eta\in [0,1]$ if for each $\nu>0$ there exists a sequence of thresholds $\{u_n^{(\nu)}\}_{n\geq 1}$ such that
\begin{equation}\label{indcruz}
nP(X_1\leq u_n^{(\nu)}<X_2)\xrightarrow [n\rain]{} \nu\qquad \textrm{and}\qquad P(\widetilde{N}_n(u_n^{(\nu)})((0,1])=0)\xrightarrow [n\rain]{}e^{-\eta \nu}.
\end{equation}
The upcrossings index $\eta$ exists if and only if there exists the extremal index $\theta$ (Ferreira (2006) \cite{fer1}) and in this case
\begin{equation}\label{rel_eta_teta}
\eta=\frac{\tau}{\nu}\ \theta.
\end{equation}
Note that for a high level $u_n$ such that $nP(X_1>u_n)=\tau'$ and $nP(X_1\leq u_n<X_2)=\nu'$ we have $\tau'\geq \nu'\geq 0$ and therefore $\eta\geq \theta.$\vspace{0.3cm}

Let us now consider the moving maxima sequence of Ferreira (2006) \cite{fer1} in the following example:

\begin{ex}
\rm{Let $\{Y_n\}_{n\geq -2}$ be a sequence of independent uniformly distributed on [0,1] variables with common distribution function $F.$ Define the 3-dependent moving maxima sequence as \linebreak $X_n=\max\{Y_n, Y_{n-2}, Y_{n-3}\},$ $n\geq 1,$ and consider $\{Z_n\}_{n\geq 1}$ a sequence of i.i.d. random variables with common distribution  $F^3$. The underlying distribution of $X_n$ is also $F^3$ and it has an extremal index $\theta=1/3$ and upcrossings index $\eta=1/2.$ Moreover, for $u_n=1-\frac{\tau}{n},$ $\tau>0,$ it holds $nP(X_1>u_n)\xrightarrow [n\rain]{} 3\tau$ and $nP(X_1\leq u_n<X_2)\xrightarrow [n\rain]{}2\tau,$ agreeing with (\ref{rel_eta_teta}). In Figure 1 we can see the sizes equal to 3 and 2, respectively for the clusters of exceedances and upcrossings of a  high level by the max-autoregressive sequence $X_n.$ Note that $X_n$ verifies condition $D^{(3)}(u_n)$ of Chernick {\it{et al.}} (1991) \cite{cher} (see Ferreira (2006) \cite{fer1}), therefore two runs of exceedances separated by one single non-exceedance are considered in the same cluster.

A shrinkage of the largest and smallest observations for the 3-dependent sequence can in general also be observed, despite the fact that we have the same model underlying both sequences.

\begin{figure}[!ht]
\begin{center}
\includegraphics[scale=0.45]{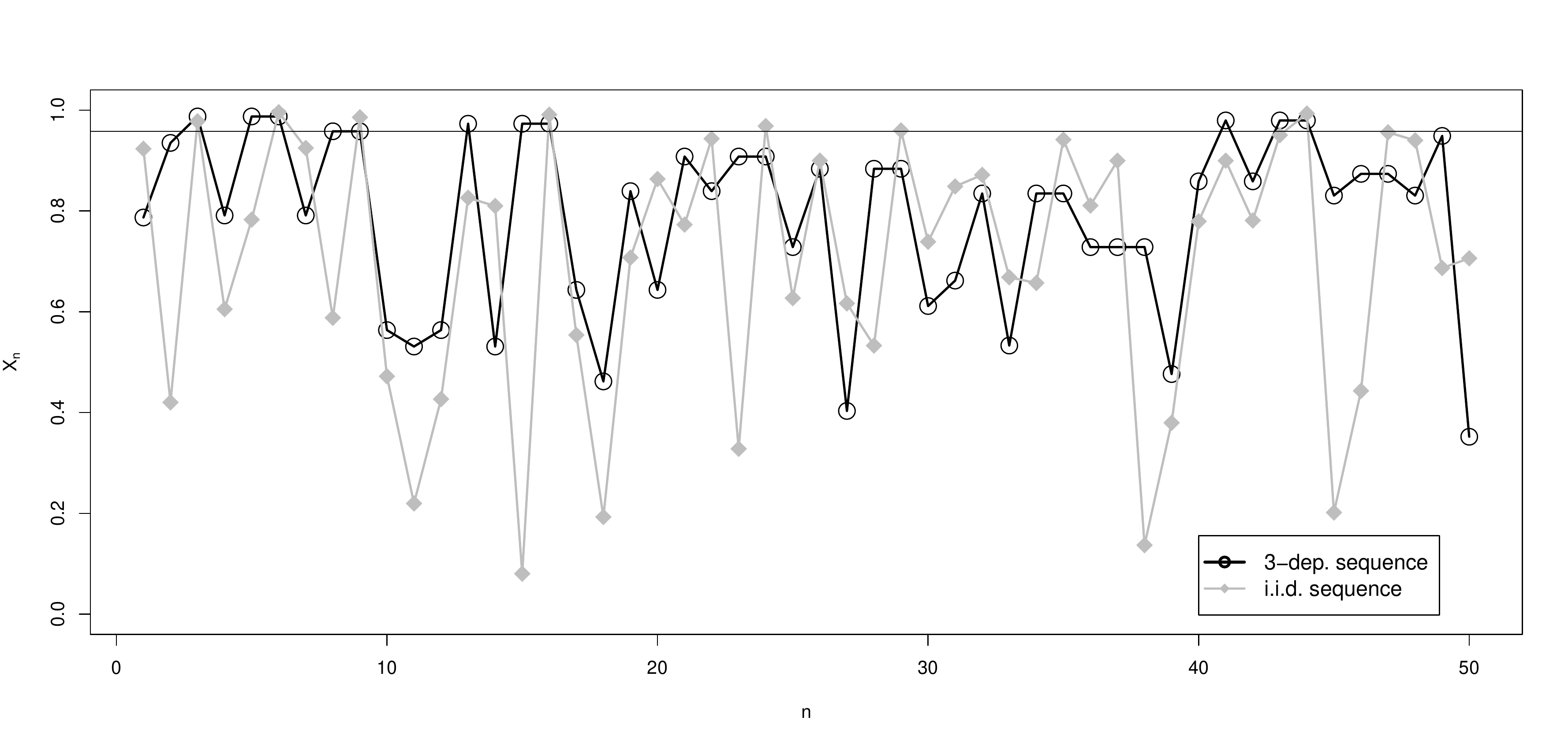}\vspace{-0.5cm}
\caption{Sample paths from an i.i.d. process and a 3-dependent process. The horizontal line corresponds the the 80th percentile of $n=50$ realizations of the 3-dependent sequence $X_n,\ n\geq 1.$ }
\end{center}
\end{figure}}
\end{ex}

In the i.i.d. setting we obviously have an upcrossings index $\eta=1.$ This also holds for processes satisfying the local dependence condition $D''(u_n)$ of Leadbetter and Nandagopalan (1989) \cite{lead2}, which enables clustering of upcrossings. Other  properties of the upcrossings index can be found in Ferreira (2006, 2007) \cite{fer1}, \cite{fer2} and Sebastião {\it{et al.}} (2010) \cite{seb}, relations with upcrossings-tail dependence coefficients are presented in Ferreira and Ferreira (2012) \cite{fer3}.

It is clear that clustering of extremes can have catastrophic consequences and upcrossings of a high threshold can also lead to adverse situations. These motivate an increasing interest in the grouping characteristics of rare events. Applications include, for example, the temporal distribution of large financial crashes or the  evaluation of oscillation periods in financial markets. A search for reliable tools to describe these features has become constant in the past years, a key aspect being the estimation of the extremal index $\theta$ and the upcrossings index $\eta$ which govern, respectively, the clustering of exceedances and upcrossings of a high level by an univariate observational series.

 Several estimators for the extremal index $\theta$ can be found in the literature (see Ancona-Navarrete and Tawn (2000) \cite{anc}, Robert {\it{et al.}} (2009) \cite{ro2} and references therein). Most of the proposed estimators are constructed by the blocks method or by the runs method. These methods identify clusters and construct estimates for $\theta$ based on these clusters. The way the clusters are identified distinguishes both methods. The blocks estimators are usually constructed by using disjoint blocks because exceedances over high thresholds of different blocks can be assumed asymptotically independent. More recently, an inter-exceedance times method has been proposed that obviates the need for a cluster identification scheme parameter (Ferro and Segers (2003) \cite{ferro}). In what concerns the estimation of the upcrossings index, little has yet been done. In Sebastião {\it{et al.}} (2013) \cite{seb1} we find
 a runs estimator and a blocks estimator derived from a disjoint blocks estimator for the extremal index, given in \cite{ro2}, for this parameter. It was shown that the runs estimator has a smaller asymptotic variance and a better performance, nevertheless it assumes the validation of a local dependence condition which does not hold for many well known processes and can be cumbersome to verify.

This paper focuses on the estimation of $\eta.$ The main novelty is our proposal of new blocks estimators and the study of their properties, namely consistency and asymptotic normality, as well as an improvement of the existing runs estimator.  We examine the behaviour of the new and existing estimators of $\eta$ and assess their performance for a range of different processes. Assuming the validation of local dependence conditions  $\widetilde{D}^{(k)}(u_n),$ $k\geq 3,$ given in Ferreira (2006) \cite{fer1} we define new runs estimators for $\eta.$ With these new estimators we obtain better estimates of $\eta,$ comparatively to the ones obtained in Sebastião {\it{et al.}} (2013) \cite{seb1} where $k=3$ was considered. We analyze a way of checking conditions $\widetilde{D}^{(k)}(u_n)$ preliminary to estimation, which helps considerably when using the runs estimators to estimate $\eta.$



We consider the problem of estimating the upcrossings index in Section 2. There we suggest an estimator motivated by the relation \begin{equation}\label{conv_eta}
E[\widetilde{N}_{r_n}(u_n^{(\nu)})\ |\ \widetilde{N}_{r_n}(u_n^{(\nu)})>0]=\sum_{j=1}^{+\infty}j
\widetilde{\pi}_n(j;u_n^{(\nu)})\xrightarrow [n\rain]{} \eta^{-1},
\end{equation}
where $$\widetilde{\pi}_n(j;u_n):=P(\widetilde{N}_{r_n}(u_n)=j\ |\ \widetilde{N}_{r_n}(u_n)>0)$$ is the conditional cluster size distribution and  $\widetilde{N}_{r_n}(u_n)\equiv\widetilde{N}_{n}(u_n)((0,r_n/n]),$ under suitable choices of $r_n\xrightarrow [n\rain]{} +\infty.$ This relation is a consequence of Lemma 2.1 of Ferreira (2006) \cite{fer1}, that states  that under condition $\Delta(u_n)$ of Hisng  {\it{et al.}} (1988) \cite{hsing} upcrossings over disjoint blocks are asymptotically independent, and definition (\ref{indcruz}) of the upcrossings index $\eta.$ It is extremely important in the interpretation of the upcrossings index, since it tells us that the reciprocal of the upcrossings index can be viewed as the limiting mean cluster size of upcrossings.

Other blocks estimators are also presented and their properties studied in Section 2. In subsection 2.2 we revisit the runs declustering method considered in Sebastião {\it{et al.}} (2013) \cite{seb1} and define new runs estimators for $\eta$ under the validation of any dependence condition $\widetilde{D}^{(k)}(u_n),$ $k\geq 3.$ A way of checking conditions $\widetilde{D}^{(k)}(u_n),$ $k\geq 3,$ preliminary to estimation is also given.

We compare, in Section 3, the performance of the estimators through a simulation study for a range of different processes with $\eta<1$ to $\eta=1$ and, in Section 4, by the application of the methods to the DAX daily log returns. Conclusions are drawn in Section 5 and proofs are given in the appendix.

\section{Estimation of the upcrossings index}

\pg We shall now consider the problem of estimating the upcrossings index in the following setting. As previously mentioned ${\bf{X}}=\{X_n\}_{n\geq 1}$ is a stationary sequence of random variables with a continuous marginal distribution function $F$ and a nonzero upcrossings index $\eta.$  Assume also that condition $\Delta(u_n^{(\nu)})$ of Hsing {\it{et al.}} (1988) \cite{hsing} holds for $\X$.

\subsection{Blocks declustering}

\pg The blocks declustering scheme consists in choosing a block length $r_n$ and partitioning the $n$ observations $(X_1,\ldots,X_n)$ of ${\bf{X}}$ into $k_n=[n/r_n]$ disjoint blocks, where $[y]$ denotes the integer part of $y.$ Each block that contains at least one upcrossing is treated as one cluster. If in Figure 1 we consider blocks of size 10 ($r_n=10$) we can identify 3 clusters.

\subsubsection{A disjoint blocks estimator}

\pg  Under relation (\ref{conv_eta}), $\eta^{-1}$ is the (asymptotic) mean cluster size of upcrossings, and so a natural way to estimate $\eta$ is as the reciprocal of the sample average cluster size, that is, count the total number of blocks with at least one upcrossing and divide it by the total number of upcrossings of a high threshold $u.$ We then propose the following blocks estimator for $\eta:$
\begin{equation}\label{estim_blocos}
\widehat{\eta}^B_n(u_n):=\frac{\sum_{i=1}^{k_n}
\indi_{\{\widetilde{N}^{(i)}_{r_n}(u_n)>0\}}}{\sum_{i=1}^{k_n}\widetilde{N}_{r_n}^{(i)}(u_n)}=
\left(\sum_{j=1}^{+\infty}
\widehat{\widetilde{\pi}}_n(j;u_n)
\right)^{-1},
\end{equation}
where
\begin{equation}
\widehat{\widetilde{\pi}}_n(j;u_n):=\sum_{i=1}^{k_n}\indi_{\{\widetilde{N}_{r_n}^
{(i)}(u_n)=j\}}\Biggl/ \sum_{i=1}^{k_n}
\indi_{\{\widetilde{N}_{r_n}^{(i)}(u_n)>0\}}\label{pi_n}
\end{equation}
is the empirical conditional distribution of $\widetilde{N}_{r_n}^{(i)}(u_n)$ and $\widetilde{N}_{r_n}^{(i)}(u_n)\equiv \widetilde{N}_{n}(u_n)([((i-1)r_n+1)/n,\ ir_n/n])=\sum_{j=(i-1)r_n+1}^{ir_n}\indi_{\{X_j\leq u_n < X_{j+1}\}},$ $1\leq i\leq k_n,$ represents the number of upcrossings of the level $u_n$ in the $ith-$block of length $r_n.$

As we can see from (\ref{estim_blocos}) and (\ref{pi_n}) the average number of upcrossings in the blocks estimates $\eta$ and the proportion of blocks with $j$ upcrossings estimates $\widetilde{\pi}(j;u_n)=\lim_{n\to +\infty}\widetilde{\pi}_n(j;u_n).$


The level $u_n$ in (\ref{estim_blocos}) is a tuning constant that determines the quality of the estimate. If $u_n$ is too high there will be just a few upcrossings of the threshold so estimates will be highly variable, but if $u_n$ is too low the approximation to the limiting characterization may be poor so the estimates may be biased. The sample size is important here as it influences what is a too high threshold.

Note that when considering $n$ finite the sequences $\{r_n\}_{n\geq 1}$ and $\{k_n\}_{n\geq 1}$ determine the cluster identification. The values of $r_n$ are determined by the dependence structure of the process. If a process has weak dependence, the high level upcrossings will occur in well defined isolated groups, so small values for $r_n$ will suffice, similarly large values of $r_n$ are required with strong dependence.

As we shall see in the next section, the quality of the estimate depends heavily on the choice of both parameters, $u_n$ and $r_n.$ 

\paragraph{Consistency and asymptotic normality}

\pg Consistency of the estimator $\widehat{\eta}^B_n$  can only be achieved  with lower thresholds than the ones considered in (\ref{indcruz}), since for these levels  there are insufficient upcrossings to give statistical ``consistency'' for the estimator. That is, as $n$ increases the value of $\widehat{\eta}^B_n$ does not necessarily converge appropriately to the value $\eta.$  We shall therefore consider thresholds  $v_n=u_{[n/c_n]}^{(\nu)},$ for some fixed $\nu>0,$ that satisfy
\begin{equation}\label{nivel}
nP(X_1\leq v_n< X_2)\sim c_n\nu \quad \textrm{as } n\to+\infty,
\end{equation}
where $\{c_n\}_{n\geq 1}$ is a sequence of real numbers such that $c_n\xrightarrow [n\rain]{}+\infty,$ $n/c_n\xrightarrow [n\rain]{}+\infty$ and
$c_n/k_n\xrightarrow [n\rain]{}0,$ for a sequence of integers  $\{k_n\}_{n\geq 1}$ satisfying $k_n\xrightarrow [n\rain]{}+\infty.$


The following result shows that for thresholds $v_n,$ satisfying (\ref{nivel}), convergence  (\ref{conv_eta}) still holds.

\begin{teo}\label{teo:blocos1}
Let $\{r_n\}_{n\geq 1},$ $\{k_n\}_{n\geq 1}$ and $\{c_n\}_{n\geq 1}$ be sequences of real numbers such that $$c_n\to +\infty,\quad k_n\to +\infty,\quad r_n\to +\infty, \quad c_n/k_n\to 0,\quad k_nr_n\leq n\quad k_nr_n\sim n,\quad {\textrm{as }} n\to +\infty$$
and $\{v_n\}_{n\geq 1}$ a sequence of thresholds defined by  (\ref{nivel}). Suppose there exists $l_n=o(r_n)$ such that $\frac{k_n}{c_n}\alpha_{n,l_n}\xrightarrow [n\rain]{}0,$ where $\alpha_{n,l_n},$ $n\geq 1,$ are the mixing coefficients of the $\Delta(v_n)$ condition.

Then
\begin{equation}
\lim_{n\to +\infty}\frac{k_n}{c_n}P(\widetilde{N}_{r_n}(v_n)>0)=\eta \nu,\label{eq:bloco1}
\end{equation} and   (\ref{conv_eta}) holds for $v_n,$ that is
$$E[\widetilde{N}_{r_n}(v_n)\ |\ \widetilde{N}_{r_n}(v_n)>0]\xrightarrow [n\rain]{} \eta^{-1}.$$
\end{teo}\vspace{0.3cm}

Note that from Theorem 2.1 if (\ref{nivel}) holds then the order of magnitude of the expected number of blocks having upcrossings of $v_n$ is $c_n\xrightarrow [n\rain]{}+\infty.$ This justifies the need to consider lower levels satisfying (\ref{nivel}) in order to guarantee the consistency of the estimator $\widehat{\eta}^B_n.$ Such levels were also considered in Sebastião {\it{et al.}} (2013) \cite{seb1} to obtain the consistency of the runs estimator of $\eta.$  \vspace{0.3cm}

Following the style of proofs used by Hsing (1991) \cite{hsing1}, we show in the next results that the blocks estimator of the upcrossings index $\widehat{\eta}^B_n$ is a consistent estimator for $\eta$ and that it is asymptotically normal. These results require assumptions on the limiting behaviour of the first and second moments of $\widetilde{N}_{r_n}(v_n).$

\begin{teo}\label{teo:consit:etaB:2} Suppose that the conditions of Theorem \ref{teo:blocos1} hold and for some $\nu>0,$
\begin{enumerate}
\item[\rm{(1)}] $\frac{k_n}{c_n}E[\widetilde{N}_{r_n}(v_n)\indi_{\{\widetilde{N}_{r_n}(v_n)>c_n\}}]\xrightarrow [n\rain]{}0,$
\item[\rm{(2)}] $\frac{k_n}{c_n^2}E[\widetilde{N}^2_{r_n}(v_n)\indi_{\{\widetilde{N}_{r_n}(v_n)\leq c_n\}}]\xrightarrow [n\rain]{}0.$
\end{enumerate}
Then
 $$\widehat{\eta}_n^B({v}_n)\xrightarrow [n\rain]{P} \eta.$$
\end{teo}\vspace{0.3cm}

The underlying idea to obtain the asymptotic distribution of the estimator $\widehat{\eta}^B_n$ is to split each block into a small block of length $l_n$ and a big block of length $r_n-l_n.$ Since $l_n=o(r_n)$ and $\frac{k_n}{c_n}\alpha_{n,l_n}\xrightarrow [n\rain]{} 0$ it is ensured that $l_n$ is sufficiently large such that blocks that are not adjacent are asymptotically independent, but does not grow to fast such that the contributions of the small blocks are negligible.






\begin{teo}\label{prop:normZ}
If conditions of Theorem \ref{teo:blocos1} hold and for some $\nu>0,$
\begin{enumerate}
\item[\rm{(1)}] $\frac{k_n}{c_n}E[\widetilde{N}^2_{r_n}(v_n)\indi_{\{\widetilde{N}^2_{r_n}(v_n)>\epsilon c_n\}}]\xrightarrow [n\rain]{}0,$ for all $\epsilon>0,$
\item[\rm{(2)}] $\frac{k_n}{c_n}E[\widetilde{N}^2_{r_n}(v_n)]\xrightarrow [n\rain]{}\eta \nu \sigma^2,$ for some $\sigma^2<+\infty.$
\end{enumerate}
Then
\begin{equation}
c_n^{-1/2}\left[\begin{tabular}{l} $\displaystyle{\sum_{i=1}^{k_n}\left(\widetilde{N}_{r_n}^{(i)}(v_n)-
E[\widetilde{N}_{r_n}^{(i)}(v_n)]\right)}$ \\
$\displaystyle{\sum_{i=1}^{k_n}\left(\indi_{\{\widetilde{N}_{r_n}^{(i)}(v_n)>0\}}-
E[\indi_{\{\widetilde{N}_{r_n}^{(i)}(v_n)>0\}}]\right)}$
\end{tabular}\right]\xrightarrow [n\rain]{d} \mathcal{N}
\left(\left[\begin{tabular}{l} $0$\\
$0$\end{tabular}\right],\left[\begin{tabular}{cc} $\nu\eta\sigma^2$ & $\nu$
\\ $\nu$ & $\nu\eta$ \end{tabular}\right]\right).\label{conv_normal}
\end{equation}
\end{teo}\vspace{0.3cm}

\begin{nota}
From assumption {\rm{(2)}} of the previous theorem we have that $E[\widetilde{N}_{r_n}^{2}(v_n)\ |\ \widetilde{N}_{r_n}^{2}(v_n)>0]\xrightarrow [n\rain]{}\sigma^2_B,$ for some $\sigma^2_B>0,$ which means that $\textrm{Var}[\widetilde{N}_{r_n}^{2}(v_n)\ |\ \widetilde{N}_{r_n}^{2}(v_n)>0]\xrightarrow [n\rain]{}\sigma^2_B-\eta^{-2},$ the asymptotic variance of a cluster size.
\end{nota}\vspace{0.3cm}

The asymptotic normality of $\widehat{\eta}^{B}_n$ is now an immediate consequence of the previous result.

\begin{coro}\label{cor:normZ}
 If the conditions of Theorem 2.3 hold then
 $$\sqrt{c_n}(\widehat{\eta}^{B}_n-\eta_n)\xrightarrow [n\rain]{d}
 \mathcal{N}\left(0,\frac{\eta}{\nu}(\eta^2\sigma^2-1)\right),$$
 where $\eta_n=E[\indi_{\{\widetilde{N}_{r_n}(v_n)>0\}}]/E[\widetilde{N}_{r_n}(v_n)].$
\end{coro}

From the previous result it is clear that the asymptotic distribution of $\widehat{\eta}^{B}_n$ is not well defined when $\sigma^2\leq 1/\eta^2,$ since in this case a smaller or equal to zero variance is obtained.\vspace{0.3cm}

\begin{nota}
The asymptotic distribution found in Corollary 2.1 differs from the one obtained with the runs estimator in Sebastião {\it{et al.}} (2013) \cite{seb1}, only on the value of $\sigma^2.$ Thus, the value of $\sigma^2$ determines which of the two is the most efficient estimator of $\eta.$
\end{nota}\vspace{0.3cm}

In order to control the bias of the upcrossings index estimator $\widehat{\eta}^{B}_n,$ lets assume that the block sizes are sufficiently large so that as $n\to +\infty$
\begin{equation}
\sqrt{c_n}(\eta_n-\eta)\xrightarrow [n\rain]{}0.\label{conv}
\end{equation}
The asymptotic variance of $\widehat{\eta}^{B}_n$ will depend on $\eta,$ $\nu$ and $\sigma^2$ which can easily be estimated with consistent estimators as shown in the next result. This result enables the construction of approximate  confidence intervals or a hypothesis test regarding $\eta.$ From a practical viewpoint it is more useful than  Corollary 2.1 and follows readily from the previous results.

\begin{coro}\label{cor:normZ2}
 If the conditions of Corollary 2.1 hold, (\ref{conv}) and $$\frac{k_n}{c_n^2}E[\widetilde{N}_{r_n}^4(v_n)\indi_{\{\widetilde{N}_{r_n}^2(v_n)\leq c_n\}}]\xrightarrow [n\rain]{}0.$$
 Then
 \begin{equation}
 \sqrt{\frac{\sum_{i=1}^{k_n}\widetilde{N}_{r_n}^{(i)}(v_n)}{\widehat{\eta}_n^B((\widehat{\eta}_n^{B}\widehat{\sigma}_n)^2-1)}}
 \ (\widehat{\eta}_n^B-\eta)\xrightarrow [n\rain]{} \mathcal{N} (0,1),\label{conv2}
 \end{equation}\vspace{0.3cm} where
 $\widetilde{N}_{r_n}^{(i)}(v_n)=\sum_{j=(i-1)r_n+1}^{ir_n}\indi_{\{X_j\leq v_n<X_{j+1}\}}$ and $\widehat{\sigma}_n^2=\frac{\sum_{i=1}^{k_n}(\widetilde{N}_{r_n}^{(i)}(v_n))^2}
 {\sum_{i=1}^{k_n}
 \indi_{\{\widetilde{N}_{r_n}^{(i)}(v_n)>0\}}}.$
 \end{coro}\vspace{0.3cm}


Given a random sample $(X_1, \ldots ,X_n)$ of ${\bf{X}}$ and a threshold $u$ we obtain, from Corollary 2.2, the asymptotic
$100(1-\alpha)$ percent two-sided confidence intervals for $\eta$
\begin{equation}
\left(\widehat{\eta}^B_n-z_{1-\alpha/2}\sqrt{\frac{\widehat{\eta}^B_n((\widehat{\eta}^B_n
\widehat{\sigma}_n)^2-1)}{\sum_{i=1}^{k_n}\widetilde{N}_{r_n}^{(i)}(u)}},\ \widehat{\eta}^B_n+z_{1-\alpha/2}\sqrt{\frac{\widehat{\eta}^B_n((\widehat{\eta}^B_n
\widehat{\sigma}_n)^2-1)}{\sum_{i=1}^{k_n}\widetilde{N}_{r_n}^{(i)}(u)}}\right),\label{int}
\end{equation}
where $z_{1-\alpha/2}$ denotes the $1-\alpha/2$ standard normal quantile, $\widehat{\eta}^B_n=\frac{\sum_{i=1}^{k_n}\indi_{\{\widetilde{N}_{r_n}^{(i)}(u)>0\}}}
{\sum_{i=1}^{k_n}\widetilde{N}_{r_n}^{(i)}(u)}$ and $\widehat{\sigma}_n^2=\frac{\sum_{i=1}^{k_n}(\widetilde{N}_{r_n}^{(i)}(u))^2}
 {\sum_{i=1}^{k_n}
 \indi_{\{\widetilde{N}_{r_n}^{(i)}(u)>0\}}}.$

The intervals (\ref{int}) are approximations of the true confidence intervals for finite samples.\vspace{0.3cm}

The finite sample properties of $\widehat{\eta}_n^B$ are now analyzed in simulated samples  from three different processes with $\eta<1$ to $\eta=1.$ The results for all examples were obtained from 5000 simulated sequences of size $n=5000,$ of each process. In each example, the tables contain Monte-Carlo approximations to biases and root mean squared errors (RMSE), in brackets, of the estimator $\widehat{\eta}_n^{B}(u)$ and the 95\% confidence intervals based on the asymptotic Normal approximation for the blocks estimates (\ref{int}), for different values of $r_n$ and different thresholds corresponding to the range between the 0.7 and 0.99 empirical quantiles, $q_{\bullet}.$ Note that the confidence intervals (\ref{int}) can only be computed when $\widehat{\sigma}^2_n>1/(\widehat{\eta}^B_n)^2.$

\begin{ex}[MM process]
{\rm{Let $X_n=\max\{Y_n,\ Y_{n-2},\ Y_{n-3}\},\ n\geq 1,$ denote the moving maxima process,
where $\{Y_n\}_{n\geq -2}$ is a sequence of independent and uniformly distributed on $[0,1]$ variables, considered in Example 1.1.

\begin{table}[!htb]
\begin{center}\scriptsize
 \begin{tabular}{lccccc}
 \hline \hline
 $r_n$& $q_{0.7}$& $q_{0.8}$& $q_{0.9}$&$q_{0.95}$&$q_{0.99}$\\ \hline
5 & 0.1983 (0.1986)  & 0.1956 (0.1960)  & 0.1964 (0.1972)  & 0.1983 (0.2000)  & 0.2001 (0.2093) \\
[-0.1cm]
 & (0.6775,0.7191) & (0.6716,0.7196) & (0.6641,0.7287) & (0.6536,0.7430) & (0.6029,0.7974)\\
7 & 0.0991 (0.0997)  & 0.1081 (0.1088)  & 0.1228 (0.1240)  & 0.1322 (0.1345)  & 0.1402 (0.1514) \\
[-0.1cm]
 & (0.5788,0.6195) & (0.5850,0.6312) & (0.5927,0.6530) & (0.5916,0.6729) & (0.5546,0.7258)\\
10 & -0.0034 (0.0102)$^+$  & 0.0215 (0.0244)  & 0.0561 (0.0583)  & 0.0770 (0.0801)  & 0.0948 (0.1069) \\
[-0.1cm]
 & (0.4765,0.5167) & (0.4984,0.5446) & (0.5266,0.5856) & (0.5389,0.6151) & (0.5208,0.6687)\\
15 & -0.1154 (0.1157)  & -0.0726 (0.0733)  & -0.0126 (0.0196)$^+$  & 0.0250 (0.0326)  & 0.0583 (0.0727)
\\ [-0.1cm]
 & (0.3658,0.4035) & (0.4042,0.4506) & (0.4570,0.5178) & (0.4873,0.5627) & (0.4960,0.6207)\\
25 & -0.2418 (0.2419)  & -0.1882 (0.1883)  & -0.0995 (0.1004)  & -0.0369 (0.0420)  & 0.0243 (0.0451)
\\ [-0.1cm]
 & (0.2430,0.2733) & (0.2906,0.3331) & (0.3688,0.4322) & (0.4229,0.5032) & (0.4699,0.5787)\\
50 & -0.3652 (0.3652)  & -0.3228 (0.3228)  & -0.2246 (0.2249)  & -0.1300 (0.1314)  & -0.0136
(0.0385) \\ [-0.1cm]
 & (0.1261,0.1436) & (0.1625,0.1919) & (0.2461,0.3047) & (0.3261,0.4138) & (0.4296,0.5432)\\
100 & -0.4324 (0.4324)  & -0.4095 (0.4095)  & -0.3403 (0.3403)  & -0.2400 (0.2404)  & -0.0568
(0.0693) \\ [-0.1cm]
 & (0.0632,0.0719) & (0.0826,0.0984) & (0.1392,0.1801) & (0.2193,0.3007) & (0.3713,0.5152)\\
1000 & -0.4932 (0.4932)  & -0.4909 (0.4909)  & -0.4836 (0.4836)  & -0.4685 (0.4685)  & -0.3476
(0.3478) \\ [-0.1cm]
 & (0.0064,0.0071) & (0.0084,0.0097) & (0.0145,0.0183) & (0.0260,0.0369) & (0.0946,0.2102)\\
 \hline \hline
\end{tabular}\caption{Bias and root mean squared error (RMSE),  in brackets, of the estimator $\widehat{\eta}_n^{B}(u)$ (first line) and the 95\% confidence intervals based on the asymptotic Normal approximation for the blocks estimates (\ref{int}) (second line), for an MM process with $\theta=1/3$ and $\eta=1/2$. Best values are ticked with a plus.}
 \end{center}
 \end{table} }}
\end{ex}

\begin{ex}[AR(1) process]
{\rm{Let  $X_n=-\frac{1}{\beta}X_{n-1}+\epsilon_n,\ n\geq 1,$ denote a negatively correlated autoregressive process of order one,
 where $\{\epsilon_n\}_{n\geq 1}$ is a sequence of i.i.d. random variables, such that, for a fixed integer $\beta \geq 2,$ $\epsilon_n\sim U\{\frac{1}{\beta},\ldots,\frac{\beta-1}{\beta},1\}$ and  $X_0\sim U(0,1)$ independent of $\epsilon_n.$

Condition $\widetilde{D}^{(3)}(u_n)$ holds for this stationary sequence, it has extremal index  $\theta=1-1/\beta^2$ and upcrossings index $\eta^{(1,{\bf{X}})}=1-1/\beta^2,$ $\beta\geq 2$ (see Sebastião {\it{et al.}} (2010) \cite{seb}).

\begin{table}[!htb]
\begin{center}\scriptsize
 \begin{tabular}{lccccc}
 \hline \hline
 $r_n$& $q_{0.7}$& $q_{0.8}$& $q_{0.9}$&$q_{0.95}$&$q_{0.99}$\\ \hline
5 & -0.1290 (0.1291)  & 0.0143 (0.0170)  & 0.0881 (0.0894)  & 0.1008 (0.1030)  & 0.1031 (0.1128) \\
[-0.1cm]
 & (0.6064,0.6355) & (0.7438,0.7848) & (0.8092,0.8670) & (0.8105,0.8910) & (0.7659,0.9402)\\
7 & -0.2834 (0.2834)  & -0.1245 (0.1248)  & 0.0098 (0.0190)$^+$  & 0.0566 (0.0611)  & 0.0747 (0.0902) \\
[-0.1cm]
 & (0.4547,0.4784) & (0.6052,0.6457) & (0.7279,0.7916) & (0.7617,0.8514) & (0.7281,0.9213)\\
10 & -0.4176 (0.4176)  & -0.2732 (0.2732)  & -0.0874 (0.0888)  & 0.0009 (0.0235)$^+$  & 0.0497 (0.0730)
\\ [-0.1cm]
 & (0.3237,0.3411) & (0.4591,0.4944) & (0.6290,0.6962) & (0.7022,0.7996) & (0.6964,0.9030)\\
15 & -0.5278 (0.5278)  & -0.4197 (0.4197)  & -0.2139 (0.2143)  & -0.0776 (0.0810)  & 0.0220 (0.0605)
\\ [-0.1cm]
 & (0.2165,0.2279) & (0.3172,0.3435) & (0.5037,0.5685) & (0.6211,0.7237) & (0.6635,0.8805)\\
25 & -0.6166 (0.6166)  & -0.5500 (0.5500)  & -0.3762 (0.3763)  & -0.1994 (0.2004)  & -0.0161
(0.0598) \\ [-0.1cm]
 & (0.1300,0.1367) & (0.1919,0.2081) & (0.3471,0.4004) & (0.4996,0.6015) & (0.6197,0.8482)\\
50 & -0.6833 (0.6833)  & -0.6500 (0.6500)  & -0.5508 (0.5508)  & -0.3889 (0.3891)  & -0.0881
(0.1041) \\ [-0.1cm]
 & (0.0650,0.0683) & (0.0960,0.1040) & (0.1834,0.2149) & (0.3197,0.4025) & (0.5414,0.7823)\\
100 & -0.7167 (0.7167)  & -0.7000 (0.7000)  & -0.6500 (0.6500)  & -0.5519 (0.5519)  & -0.1983
(0.2041) \\ [-0.1cm]
 & (0.0325,0.0341) & (0.0480,0.0520) & (0.0921,0.1080) & (0.1727,0.2234) & (0.4311,0.6723)\\
1000 & -0.7467 (0.7467)  & -0.7450 (0.7450)  & -0.7400 (0.7400)  & -0.7300 (0.7300)  & -0.6500
(0.6500) \\ [-0.1cm]
 & (0.0033,0.0034) & (0.0048,0.0052) & (0.0093,0.0107) & (0.0178,0.0223) & (0.0712,0.1288)\\
 \hline \hline
\end{tabular}\caption{Bias and root mean squared error (RMSE),  in brackets, of  $\widehat{\eta}_n^{B}(u)$ (first line) and the 95\% confidence intervals based on the asymptotic Normal approximation for the blocks estimates (\ref{int}) (second line), for a negatively correlated AR(1) process with $\beta=2$ and  $\theta=\eta=0.75$. Best values are ticked with a plus.}
 \end{center}
 \end{table}}}
\end{ex}

\begin{ex}[MAR(1) process]
{\rm{Let $X_n=\alpha \max\{X_{n-1},\ \epsilon_n\},\ n\geq 1,\ 0<\alpha<1,$ denote the max-autoregressive process of order one where $X_0$ is a random variable  with d.f.
$H_0,$ independent of the sequence $\{\epsilon_n\}_{n\geq 1}$ of i.i.d. random variables with unit Fréchet distribution. This process, which is a special case of the general MARMA(p,q) processes introduced by Davis and Resnick (1989) \cite{davis}, verifies condition  $D''(u_n),$ has extremal index $\theta=1-\alpha,\ 0<\alpha<1$ (see Alpuim (1989) \cite{al} and Ferreira (1994) \cite{fer0}) and upcrossings index $\eta=1.$

\begin{table}[!htb]
\begin{center}\scriptsize
 \begin{tabular}{lccccc}
  \hline \hline
 $r_n$& $q_{0.7}$& $q_{0.8}$& $q_{0.9}$&$q_{0.95}$&$q_{0.99}$\\ \hline
5 & -0.0061 (0.0092)  & -0.0039 (0.0075)  & -0.0019 (0.0065)  & -0.0008 (0.0056)$^+$  & -0.0003 (0.0075)$^+$
\\ [-0.1cm]
 & (0.9843,1.0035) & (0.9893,1.0029) & (0.9947,1.0016) & (0.9978,1.0007) & (0.9992,1.0002)\\
7 & -0.0143 (0.0177)  & -0.0096 (0.0140)  & -0.0048 (0.0110)  & -0.0024 (0.0100)  & -0.0009 (0.0122)
\\ [-0.1cm]
 & (0.9678,1.0036) & (0.9759,1.0048) & (0.9871,1.0033) & (0.9934,1.0017) & (0.9982,1.0001)\\
10 & -0.0305 (0.0341)  & -0.0202 (0.0250)  & -0.0100 (0.0174)  & -0.0050 (0.0148)  & -0.0013
(0.0149) \\ [-0.1cm]
 & (0.9413,0.9976) & (0.9546,1.0049) & (0.9736,1.0064) & (0.9861,1.0039) & (0.9968,1.0006)\\
15 & -0.0659 (0.0695)  & -0.0446 (0.0494)  & -0.0232 (0.0317)  & -0.0112 (0.0236)  & -0.0029
(0.0214) \\ [-0.1cm]
 & (0.8935,0.9747) & (0.9161,0.9948) & (0.9446,1.0090) & (0.9702,1.0074) & (0.9931,1.0012)\\
25 & -0.1495 (0.1525)  & -0.1044 (0.1090)  & -0.0549 (0.0636)  & -0.0280 (0.0427)  & -0.0055
(0.0301) \\ [-0.1cm]
 & (0.7946,0.9063) & (0.8378,0.9534) & (0.8887,1.0015) & (0.9303,1.0138) & (0.9866,1.0025)\\
50 & -0.3521 (0.3539)  & -0.2624 (0.2658)  & -0.1466 (0.1542)  & -0.0766 (0.0924)  & -0.0153
(0.0515) \\ [-0.1cm]
 & (0.5837,0.7122) & (0.6621,0.8130) & (0.7653,0.9416) & (0.8385,1.0082) & (0.9629,1.0065)\\
100 & -0.6059 (0.6066)  & -0.4966 (0.4982)  & -0.3123 (0.3177)  & -0.1751 (0.1893)  & -0.0365
(0.0838) \\ [-0.1cm]
 & (0.3447,0.4435) & (0.4328,0.5740) & (0.5818,0.7936) & (0.6988,0.9510) & (0.9148,1.0122)\\
1000 & -0.9591 (0.9592)  & -0.9431 (0.9431)  & -0.8932 (0.8933)  & -0.7921 (0.7932)  & -0.3382
(0.3843) \\ [-0.1cm]
 & (0.0363,0.0454) & (0.0487,0.0651) & (0.0835,0.1301) & (0.1420,0.2738) & (0.4526,0.8711)\\
 \hline \hline
 \end{tabular}\caption{Bias and root mean squared error (RMSE),  in brackets, of  $\widehat{\eta}_n^{B}(u)$ (first line) and the 95\% confidence intervals based on the asymptotic Normal approximation for the blocks estimates (\ref{int}) (second line), for a MAR(1) process with $\alpha=0.9,$ $\theta=0.1$ and $\eta=1.$ Best values are ticked with a plus.}
 \end{center}
 \end{table}}}
\end{ex}\vspace{0.1cm}

Absolute biases tend to be smallest at small block sizes, {\it{i.e.}}, $r_n=5,$ 7, 10 and $15.$ Furthermore, the results suggest that generally speaking the absolute bias of $\widehat{\eta}_n^{B}(u)$ increases with $r_n$ but decreases with $u$ and the variance of $\widehat{\eta}_n^{B}(u)$ increases with $u$ but decreases with $r_n.$ The quality of the estimate depends strongly on the choice of the two parameters $u$ and $r_n,$ as previously stated. Nevertheless, the blocks estimator has a better performance at the 90\% and the 95\% thresholds. The poor performance at the 99\% threshold is justified by the few observed threshold upcrossings.

\begin{nota}
We point out from the previous examples that the stronger (weaker) the dependence between upcrossings, {\it{i.e.}} the bigger (smaller) clustering of upcrossings, which implies the smaller (bigger) $\eta$, the bigger (smaller) values of $r_n$ are required.
\end{nota}

\begin{nota}
In the estimation of the extremal index the block size $r_n$ has commonly been chosen as the square root of the sample size, {\it{i.e.}} $r_n\approx \sqrt{n}$ (Gomes (1993) \cite{gom0}, Ancona-Navarrete and Tawn (2000) \cite{anc}, among others). The previous examples show that such a large block size is not a reasonable choice in the estimation of the upcrossings index. Choices ranging from $r_n=[n/(\log n)^{3.2}]$ to $r_n=[n/(\log n)^{2.7}]$ seem more reasonable.
\end{nota}


\paragraph{Some considerations about the choice of the levels}

\pg In practical applications the deterministic levels $v_n$ previously considered will have to be estimated, {\it{i.e.}}, they will have to be replaced by random levels suggested by the relation $\frac{n}{c_n}P(X_1\leq v_n<X_2)\sim \nu,$ as $n\to +\infty.$ Nevertheless, these random levels cannot be represented by an appropriate order statistics, contrarily to the random levels used in the estimation of the extremal index.

As in Sebastião {\it{et al.}} (2013) \cite{seb1} we shall consider the random levels $\widehat{v}_n=X_{n-[c_n\tau]:n},$ used in the estimation of the extremal index, more precisely, we shall consider  the blocks estimator
 \begin{equation}\label{estim_blocos_al}
\widetilde{\eta}^B_n(\widehat{v}_n):=\frac{\sum_{i=1}^{k_n}
\indi_{\{\widetilde{N}^{(i)}_{r_n}(\widehat{v}_n)>0\}}}{\sum_{i=1}^{k_n}\widetilde{N}_{r_n}^{(i)}(\widehat{v}_n)}=
\left(\sum_{j=1}^{+\infty}j
\widehat{\widetilde{\pi}}_n(j;\widehat{v}_n)
\right)^{-1}.
\end{equation}
For these random levels $\widehat{v}_n,$ we show in the following result that the blocks estimator $\widetilde{\eta}^B_n(\widehat{v}_n)$ is also a consistent estimator of $\eta.$

\begin{teo}\label{teo:aleatorio}
Suppose that for each $\nu>0$ there exists $v_n=u_{n/c_n}^{(\nu)}=u_{n/c_n}^{(\tau)}$ ($\frac{n}{c_n}P(X_1\leq v_n <X_2)\sim \nu$ and $\frac{n}{c_n}P(X_1>v_n)\sim \tau$ as $n\to +\infty$) for some $\tau>0,$ the conditions of Theorem 2.2 hold for all $\nu>0,$ then
$$\widehat{\eta}^{B}_n(\widehat{v}_n)\xrightarrow [n\rain]{}\eta.$$
\end{teo}\vspace{0.3cm}

\begin{nota}
If ${\bf{X}}$ has extremal index $\theta>0$ and upcrossings index $\eta>0,$ then there exists $v_n=u_{n/c_n}^{(\nu)}=u_{n/c_n}^{(\tau)}$ and $\eta\nu=\theta\tau$ (Ferreira (2006) \cite{fer1}).
\end{nota}\vspace{0.3cm}

The high level $v_n,$ considered in Theorem 2.4, must be such that $nP(X_1>v_n)=c_n\tau=\tau_n,$ $\tau_n\xrightarrow [n\rain]{}+\infty$ and $\tau_n/n\xrightarrow [n\rain]{}0.$ In the simulation study, presented further ahead, we shall consider the deterministic levels $u_n=X_{n-s:n},$ corresponding to the ($s+1$)th top order statistics associated to a random sample $(X_1,\ldots,X_n)$ of ${\bf{X}}.$ The upcrossings index estimator will become in this way a function o $s$. In fact, $s$ is replacing $\tau_n=c_n\tau,$ $c_n\xrightarrow [n\rain]{}+\infty.$ Consistency is attained only if $s$ is intermediate, {\it{i.e.}}, $s\to+\infty$ and $s=o(n)$ as $n\to +\infty.$

The choice of $s,$ the number of top order statistics to be used in the estimation procedure is a sensitive and complex problem in extreme value applications. In Neves {\it{et al.}} (2014) (\cite{neves}) an heuristic algorithm is used for the choice of the optimal sample fraction in the estimation of the extremal index. The optimal sample fraction is chosen by looking at the sample path of the estimator of the extremal index as a function of $s$ and finding the largest range of $s$ values associated with stable estimates. Since the estimates of the extremal index and the upcrossings index sometimes do not have the largest stability region around the target value, this algorithm may give misleading results.

\subsubsection{Short note on other blocks estimators}

\pg If the stationary sequence $\X$ also has extremal index $\theta>0$ and there exists a sequence of thresholds $\{u_{r_n}\}_{r_n\geq 1}$ such that $\displaystyle\lim_{r_n\to +\infty}r_n P(X_1>u_{r_n})=\tau$ and $\displaystyle\lim_{r_n\to +\infty}r_n P(X_1\le u_{r_n}<X_2)=\nu$ then, since (\ref{rel_eta_teta}) holds, it follows that
\begin{equation}
\eta_{r_n}(u_{r_n})=-\frac{\log F_{r_n}(u_{r_n})}{r_nP(X_1\le u_{r_n}<X_2)} \xrightarrow [r_n\rain]{} \eta,\label{rel_teta}
\end{equation}
where $M_{r_n}=\max\{X_1,\ldots,X_{r_n}\}$ and $F_{r_n}(u_{r_n})=P(M_{r_n}\leq u_{r_n}).$

The distribution function of the block maximum $F_{r_n}$ in (\ref{rel_teta}) can be estimated using maxima of $k_n=[n/r_n]$ disjoint blocks or using maxima of $n-r_n+1$ sliding blocks, as suggested by Robert {\it{et al.}} (2009) \cite{ro2}, in the following ways
\begin{equation*}
 \widehat{F}^{dj}_{n}(u_{r_n})=\frac{1}{k_n}\sum_{i=1}^{k_n}\indi_{\{M_{(i-1)r_n+1,ir_n}\le u_{r_n}\}},\quad  \widehat{F}^{sl}_{n}(u_{r_n})=\frac{1}{n-r_n+1}\sum_{i=1}^{n-r_n+1}\indi_{\{M_{i,r_n+i-1}\le u_{r_n}\}}
\end{equation*}
where $M_{i,j}=\max\{X_i,\ldots,X_j\},$ for $0<i<j.$ Robert {\it{et al.}} (2009) \cite{ro2} showed that $\widehat{F}_{n}^{dj}(u_{r_n})-\widehat{F}_{n}^{sl}(u_{r_n})$ has a non-negligible asymptotic variance and is asymptotically uncorrelated with $\widehat{F}_{n}^{sl}(u_{r_n}).$ Thus, the sliding blocks estimator is the most efficient convex combination of the disjoint and sliding blocks estimators for $F_{r_n}(u_{r_n}).$

We can now define a disjoint and a sliding blocks estimator for the upcrossings index  $\eta,$ respectively, as
\begin{equation}\label{estim_dj}
\widehat{\eta}^{dj}_{n}(u_{r_n}):=-\frac{\log \widehat{F}_{n}^{dj}(u_{r_n})}{\frac{1}{k_n}\sum_{i=1}^{k_n}\widetilde{N}_{r_n}^{(i)}(u_{r_n})},
\end{equation}
and
\begin{equation}\label{estim_sl}
\widehat{\eta}^{sl}_{n}(u_{r_n}):=-\frac{\log \widehat{F}_{n}^{sl}(u_{r_n})}{\frac{1}{k_n}\sum_{i=1}^{k_n}\widetilde{N}_{r_n}^{(i)}(u_{r_n})}.
\end{equation}\vspace{0.3cm}

The estimators $\widehat{\eta}^{dj}_{n}$ and $\widehat{\eta}^{sl}_{n}$ in (\ref{estim_dj}) and (\ref{estim_sl}) are not defined when all blocks have at least one upcrossing and $\widehat{\eta}^{dj}_{n}$ has already appeared in Sebastião {\it{et al.}} (2013) \cite{seb1}.\vspace{0.3cm}

\begin{nota}
If for each $\nu>0$ there exists $u_n^{(\nu)}$ and $u_n^{(\nu)}=u_n^{(\tau)}$ for some $\tau>0,$ then $P(\max\{X_1,\ldots,X_n\}\leq u_n^{(\tau)})-P(\widetilde{N}_n(u_n^{(\nu)})=0)\xrightarrow [n\rain]{}0$ (Ferreira (2006) \cite{fer1}). Therefore, for large $n$ we can replace the probability in the numerator of (\ref{rel_teta}) with $P(\widetilde{N}_{r_n}(u_{r_n})=0)$ and estimate it by its empirical counterpart. Nevertheless, estimator (\ref{estim_dj}) creates synergies in the estimation of the pair $(\theta,\ \eta).$ We recall that the extremal index $\theta$ and the upcrossings index $\eta$ provide different but complementary information concerning the occurrence of rare events.
\end{nota}\vspace{0.3cm}

Imposing on ${\bf{X}}$ a stronger condition than condition $\Delta(u_n),$ involving its maximal correlation coefficients, and considering the threshold sequence $\{u_{r_n}\}_{r_n\geq 1}$ deterministic, we have, {\it{mutatis mutandis}}, from the arguments used in Robert {\it{et al.}} (2009) \cite{ro2}, that $\widehat{\eta}^{dj}_{n}$ and $\widehat{\eta}^{sl}_{n}$ are consistent estimators of $\eta.$ If in addition, there exists a constant $p>1$ such that $E\left[\widetilde{N}_{r_n}^{2p}(u_{r_n})\right]=O(1)$, as $n\to +\infty,$ and  $\sqrt{k_n}(\eta_{r_n}(u_{r_n})-\eta)\xrightarrow [r_n\rain]{} 0,$ then from the central limit theorem for triangular arrays it can be shown that
\begin{equation*}
 \sqrt{k_n}\left(\begin{array}{l}
 \widehat{\eta}_{n}^{dj}(u_{r_n})-\eta\\[0.1cm]
 \widehat{\eta}_{n}^{sl}(u_{r_n})-\eta
                 \end{array}\right)\xrightarrow [r_n\rain]{d} N(0,\mathbf{V})
\end{equation*}
where $\mathbf{V}=[v_{i,j}]_{2\times 2}$ is a symmetric matrix with
\begin{eqnarray*}
 v_{11}&=&\frac{\eta}{\nu}\left(\frac{e^{\eta\nu}}{\eta\nu}-\frac{1}{\eta\nu}+c^2-1\right),\\[0.1cm]
 v_{22}=v_{12}&=&\frac{\eta}{\nu}\left(2\times \frac{e^{\eta\nu}-1}{\eta^2\nu^2}-\frac{1}{\eta\nu}+c^2-1\right)
\end{eqnarray*}
and $c^2=\frac{\sum_{j\ge 1} j^2\widetilde{\pi}(j;u_{r_n})-\eta^{-2}}{\eta^{-2}}.$\vspace{0.3cm}

\begin{nota}
Since $v_{22}\leq v_{11}$ the sliding blocks estimator $\widehat{\eta}^{sl}_{n}$ is more efficient than its disjoint version. On the other hand, considering that $c_n/k_n\to 0,$ as $n\to +\infty,$ the blocks estimator $\widehat{\eta}^{B}_{n}$ is a more efficient estimator than $\widehat{\eta}_{n}^{dj}.$
\end{nota}

\subsection{Runs declsutering}

\pg Runs declustering assumes that upcrossings belong to the same cluster if they are separated by fewer than a certain number of non-upcrossings of the threshold.  More precisely, if the process $\X$ verifies condition $\widetilde{D}^{(k)}(u_n),$ for some $k\geq 2,$ of Ferreira (2006) \cite{fer1}, that locally restricts the dependence of the sequence but still allows clustering of upcrossings, runs of upcrossings in the same cluster must be separated at most by $k-2$ non-upcrossings. In Figure 1 we can identify 3 clusters of size 2, since the process verifies condition $\widetilde{D}^{(3)}(u_n).$\vspace{0.3cm}

{\bf{Condition $\widetilde{D}^{(k)}(u_n),$ $k\geq 2,$}} is said to be satisfied if
\begin{equation}
\lim_{n\to
+\infty}nP(X_1\leq u_n<X_2,\ \widetilde{N}_{3,k}=0,\
\widetilde{N}_{k+1,r_n}>0)=0,\label{Dk}
\end{equation}
for some sequence $r_n=[n/k_n]$ with
$\{k_n\}_{n\geq 1}$ satisfying $k_n\xrightarrow [n\rain]{}+\infty,$ $\frac{k_nl_n}{n}
\xrightarrow [n\rain]{}0,$ $k_n\alpha_{n,l_n}\xrightarrow
[n\rain]{} 0,$ where $\alpha_{n,l_n}$ are the mixing coefficients of the $\Delta(u_n)$ condition, and $\widetilde{N}_{i,j}\equiv\widetilde{N}_n(u_n)([i/n,j/n])$ with $\widetilde{N}_{i,j}=0$ if $j<i.$

This family of local dependence conditions is slightly stronger than $D^{(k)}(u_n)$ of Chernick {\it{et al.}} (1991) \cite{cher} and  (\ref{Dk}) is implied by
$$n\sum_{i=k+1}^{r_n-1}P(X_1\leq u<X_2,\ X_i\leq u<X_{i+1})\xrightarrow
[n\rain]{} 0. $$
When $k=2$ we find the slightly weakened condition  $D''(u_n)$ of Leadbetter and Nandagopalan (1989) \cite{lead2}.

Under this hierarchy of increasingly weaker mixing conditions $\widetilde{D}^{(k)}(u_n),$ $k\geq 2,$ the upcrossings index can be computed as follows:

\begin{prop}[Ferreira (2006) \cite{fer1}]
If ${\bf{X}}$ satisfies condition $\Delta(u_n)$ and for some $k\geq 2,$ condition $\widetilde{D}^{(k)}(u_n^{(\nu)})$ holds for some $\nu>0,$ then the upcrossings index of ${\bf{X}}$ exists and is equal to $\eta$ if and only if $$P(\widetilde{N}_{3,k}(u_n^{(\nu)})=0\ |\ X_1\leq u_n^{(\nu)}<X_2)\xrightarrow [n\rain]{} \eta $$ for each $\nu>0.$
\end{prop}

\subsubsection{The runs estimators}

\pg The estimators considered up to now treat each block as one cluster, but under a local dependence condition $\widetilde{D}^{(k)}(u_n),$ $k\geq 3,$ clusters can be identified in a different way, for example, runs of upcrossings or runs of upcrossings separated by at most one non-upcrossings of a certain threshold may define a cluster. More precisely, if the process $\X$ verifies condition $\widetilde{D}^{(3)}(u_n),$  upcrossing clusters may be simply identified asymptotically as runs of consecutive upcrossings and the cluster sizes as run lengths (Ferreira (2007) \cite{fer2}). Therefore, from Propositon 2.1 if for some $k\geq 3$ condition $\widetilde{D}^{(k)}(u_n)$ holds, the upcrossings index can be estimated by the ratio between the total number of $k-2$ non-upcrossings followed by an upcrossing and the total number of upcrossings, {\it{i.e.}}, the runs estimator  given by
\begin{equation}\label{estim_run}
\widehat{\eta}^{R|k}_n=\widehat{\eta}^{R|k}_n(u_n):=\frac{\sum_{i=1}^{n-k}
\indi_{\{\widetilde{N}_{i,i+k-3}=0,\ X_{i+k-1}\leq u_n<X_{i+k}\}}}{\sum_{i=1}^{n-1}\indi_{\{X_{i}\leq u_n<X_{i+1}\}}},
\end{equation}
with $\widetilde{N}_{i,j}\equiv\widetilde{N}_n(u_n)([i/n,j/n])$ and $\widetilde{N}_{i,j}=0$ if $j<i.$

When $k=3$ the runs estimator corresponds to the one fairly studied in Sebastião {\it{et al.}} (2013) \cite{seb1}. There it was shown that the runs estimator has smaller bias and mean squared error when compared to the disjoint blocks estimator $\widehat{\eta}^{dj}_n.$

Properties such as consistency and asymptotic normality can be proved for the estimators $\widehat{\eta}^{R|k}_n,$ $k\geq 4,$ with similar arguments to the ones used in Sebastião {\it{et al.}} (2013) \cite{seb1}.  Nevertheless, their validation is restricted to the validation of a local dependence condition $\widetilde{D}^{(k)}(u_n),$ $k\geq 4.$

In what follows we propose a way of checking conditions $\widetilde{D}^{(k)}(u_n),$ $k\geq 3,$ preliminary to estimation, which turns out to be the only possible solution when dealing with real data. Even when the underlying model is known, the verification of conditions $\widetilde{D}^{(k)}(u_n),$ $k\geq 3,$ can be cumbersome, so the following procedure can in these situations be an important auxiliary tool. A similar approach has been followed by Süveges (2007) \cite{suv} to check condition $D^{(2)}(u_n)$ of Chernick {\it{et al.}} (1991) \cite{cher} which is slightly stronger than condition $D''(u_n)$ of Leadbetter and Nandagopalan (1989) \cite{lead2}.\vspace{0.3cm}

Condition $\widetilde{D}^{(k)}(u_n),$ $k\geq 3,$ essential for the validity of the runs estimator $\widehat{\eta}^{R|k}_n,$ may be checked by calculating
\begin{equation}
p^{(k)}(u,r)=\frac{\sum_{i=1}^{n-r+1}\indi_{\{X_i\leq u<X_{i+1},\ \widetilde{N}_{i+2,i+k-1}=0,\ \widetilde{N}_{i+k,i+r-1}>0\}}}{\sum_{i=1}^{n-1} \indi_{\{X_i\leq u<X_{i+1}\}}}\label{proporcao}
\end{equation}
for a sample $(X_1,\ldots,X_n)$ of ${\bf{X}},$ $u$ a high threshold and $r$ the size of the blocks into which the sample is partitioned.  \vspace{0.3cm}

The limit condition $\widetilde{D}^{(k)}(u_n),$ for some $k\geq 3,$ will be satisfied if  there exists a path $(u_i,r_j)$ with $u_i\to +\infty$ and $r_j\to +\infty$ for which the $p^{(k)}(u_i,r_j)\to 0.$ With (\ref{proporcao}) we are looking for the so called anti-$\widetilde{D}^{(k)}(u_n),$ $k\geq 3,$ events $\{\widetilde{N}_{i+2,i+2}=0,\ \widetilde{N}_{i+3,i+r-1}>0\ |\ X_i\leq u<X_{i+1}\}$ among the upcrossings for a range of thresholds and block sizes. This proportion of anti-$\widetilde{D}^{(k)}(u_n),$ $k\geq 3,$ events gives us information on how many clusters are misidentified as two or more clusters and therefore information of the bias. If such clusters are few, it is plausible to accept the relatively small upward bias despite the possible failure of $\widetilde{D}^{(k)}(u_n),$ $k\geq 3,$ and use $\widehat{\eta}^{R|k}_n$ to estimate $\eta.$ This gives us a way to choose, in practical applications, among the estimators  $\widehat{\eta}^{R|k}_n,$ $k\geq 3,$ the one to estimate $\eta.$

In Figures 2 and 3 we plot the proportions $p(u,r)\equiv p^{(3)}(u,r),$  considering $u$ the empirical quantiles with probabilities ranging from 0.95 to 0.995 and $r$ ranging from 5 to 20, for sequences of size $n=10000,$ of the processes considered in the examples of Section 2, of an AR(2) process and a GARCH(1,1) process with Student-$t$ distributed innovations. For the sake of clarity we consider $\widehat{F}(u)$ instead of $u$ in the plots.

Similar results were obtained for other sample sizes, namely $n=5000$ and 30000. For the AR(2) process the proportions $p^{(4)}(u,r)$ and $p^{(5)}(u,r)$ did not differ much from $p^{(3)}(u,r)$ which suggest that this process also does not verify conditions $\widetilde{D}^{(4)}(u_n)$ and $\widetilde{D}^{(5)}(u_n).$ The anti-$\widetilde{D}^{(k)}(u_n),$ $k=3,4,5,6,$  proportions for the GARCH(1,1) process, plotted in Figure 3, show a more prominent decrease, towards zero, for $k=5$ and $k=6.$ This suggests that conditions $\widetilde{D}^{(3)}(u_n)$ and $\widetilde{D}^{(4)}(u_n)$ are unlikely to hold for this process and that $k=5$ seems a reasonable choice.

\begin{figure}[!ht]
\begin{center}
\includegraphics[scale=0.55]{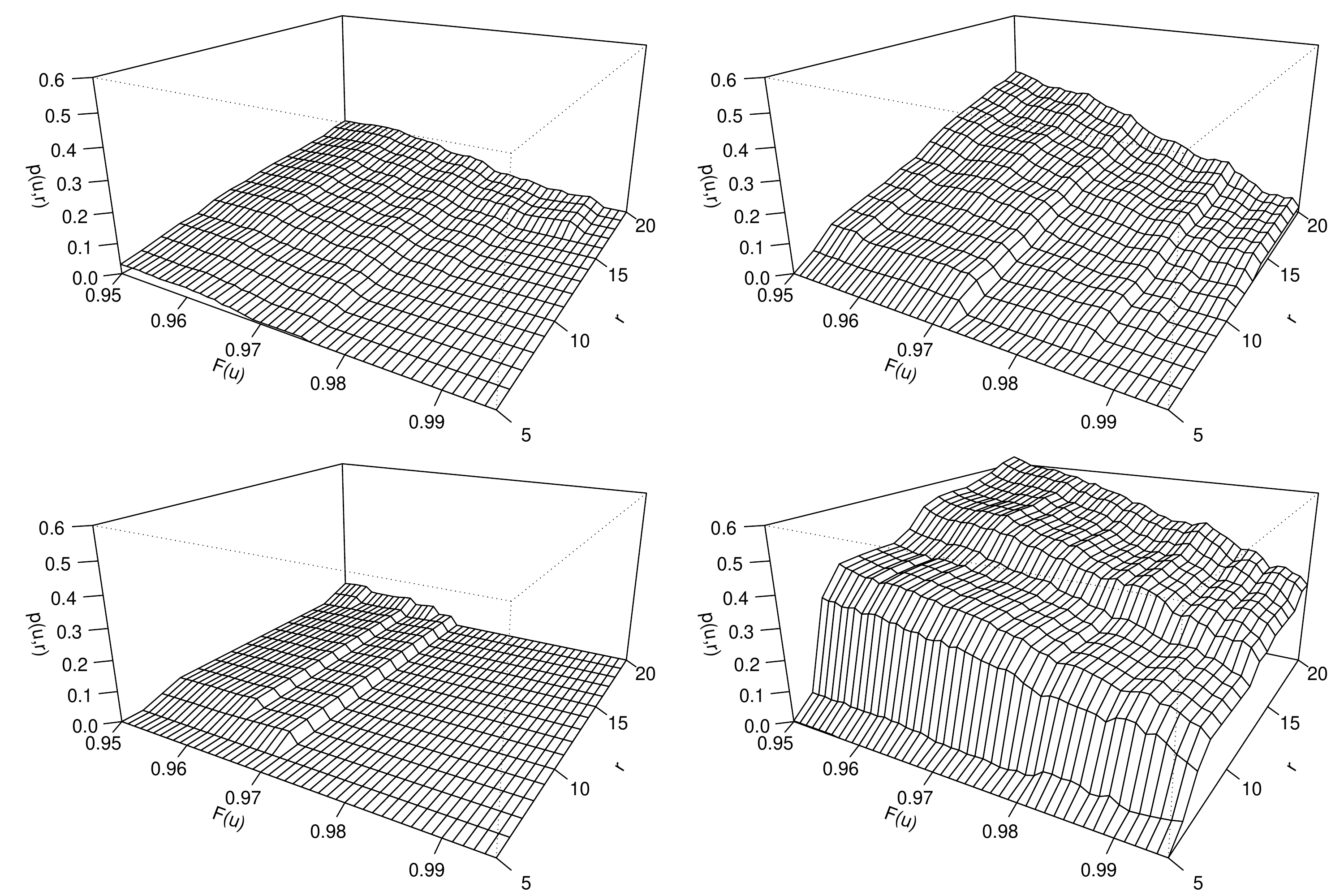}\vspace{-0.5cm}
\caption{Observed proportions $p(u,r)\equiv p^{(3)}(u,r)$ for the MM process (top left panel), negatively correlated AR(1) process with $s=2$ (top right panel), MAR(1) process (bottom left panel) and AR(2) process with $\phi_1=0.93$ and $\phi_2=-0.86$ (bottom right panel).}
\end{center}
\end{figure}

\begin{figure}[!ht]
\begin{center}
\includegraphics[scale=0.55]{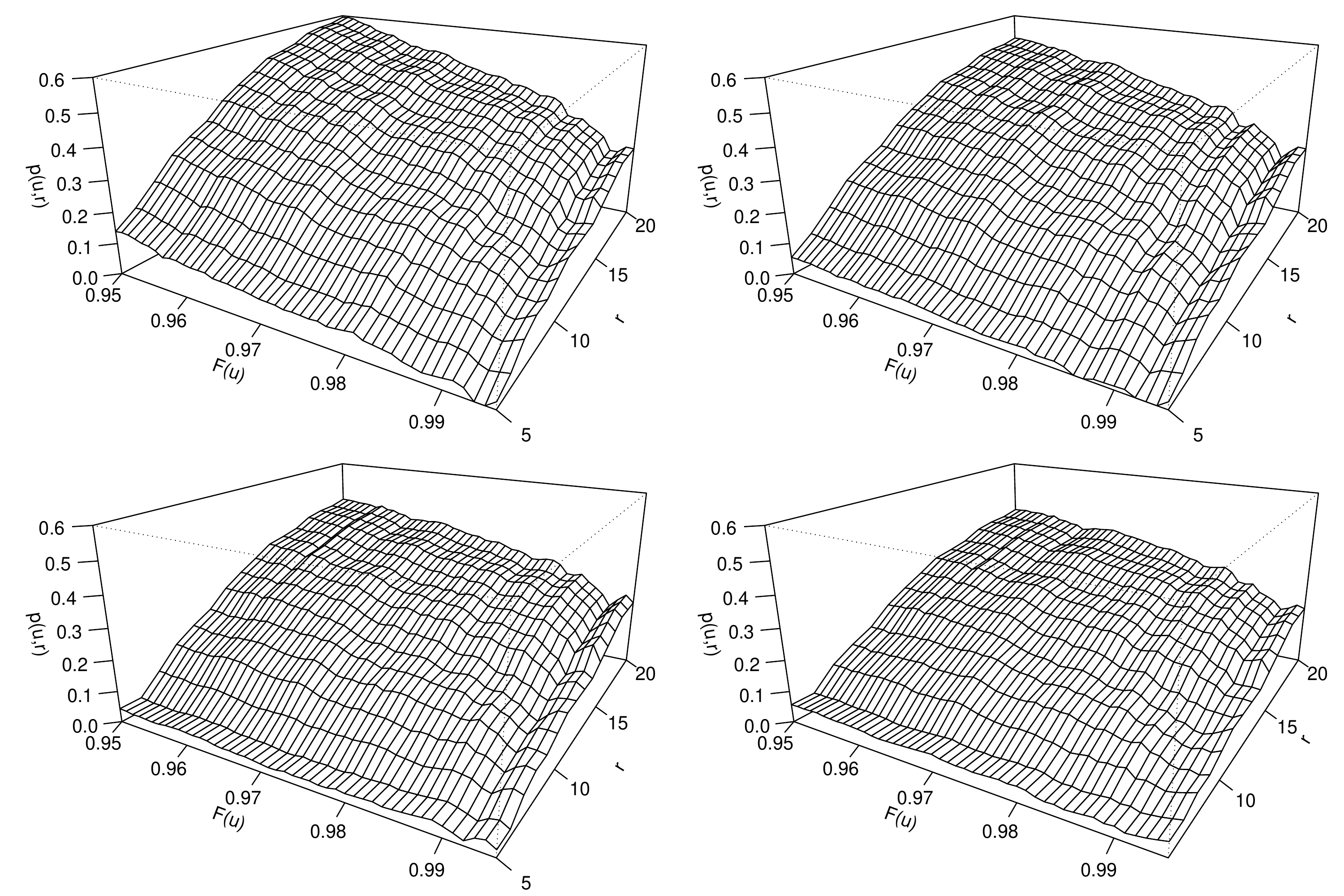}\vspace{-0.5cm}
\caption{ From left to right and top to bottom, observed proportions $p(u,r)\equiv p^{(k)}(u,r),$ $k=3,4,5$ and 6 for a  GARCH(1,1) process with autoregressive parameter $\alpha=0.08$, variance parameter $\beta = 0.87$ and innovations $t_7$}
\end{center}
\end{figure}

\vspace{0.3cm}We recall that the conditions $\widetilde{D}^{(k)}(u_n),$ $k\geq 3,$ form a hierarchy of increasingly weaker mixing conditions and therefore, if a process verifies condition $\widetilde{D}^{(k')}(u_n)$ for a fixed $k',$ then it verifies all conditions $\widetilde{D}^{(k)}(u_n)$ with $k\geq k'.$ Accordingly we can use any of the estimators $\widehat{\eta}^{R|k}_n,$ $k=k',k'+1,\ldots,$ to estimate $\eta.$ We simulated for the estimators $\widehat{\eta}^{R|k}_n,$ $k=3,4,5,$ the mean (E) and the root mean squared error (RMSE) from 5000 samples of sizes $n=5000$ of the MM process given in Example 2.1. In Figure 4, we find the estimated values as functions of $s\geq 1,$ corresponding to the estimates obtained with $\widehat{\eta}^{R|k}_n(X_{n-s:n}),$ $k=3,4,5,$ where $X_{n-s:n}$ denotes the ($s+1$)th top order statistics.

\begin{figure}[!ht]
\begin{center}
\includegraphics[scale=0.35]{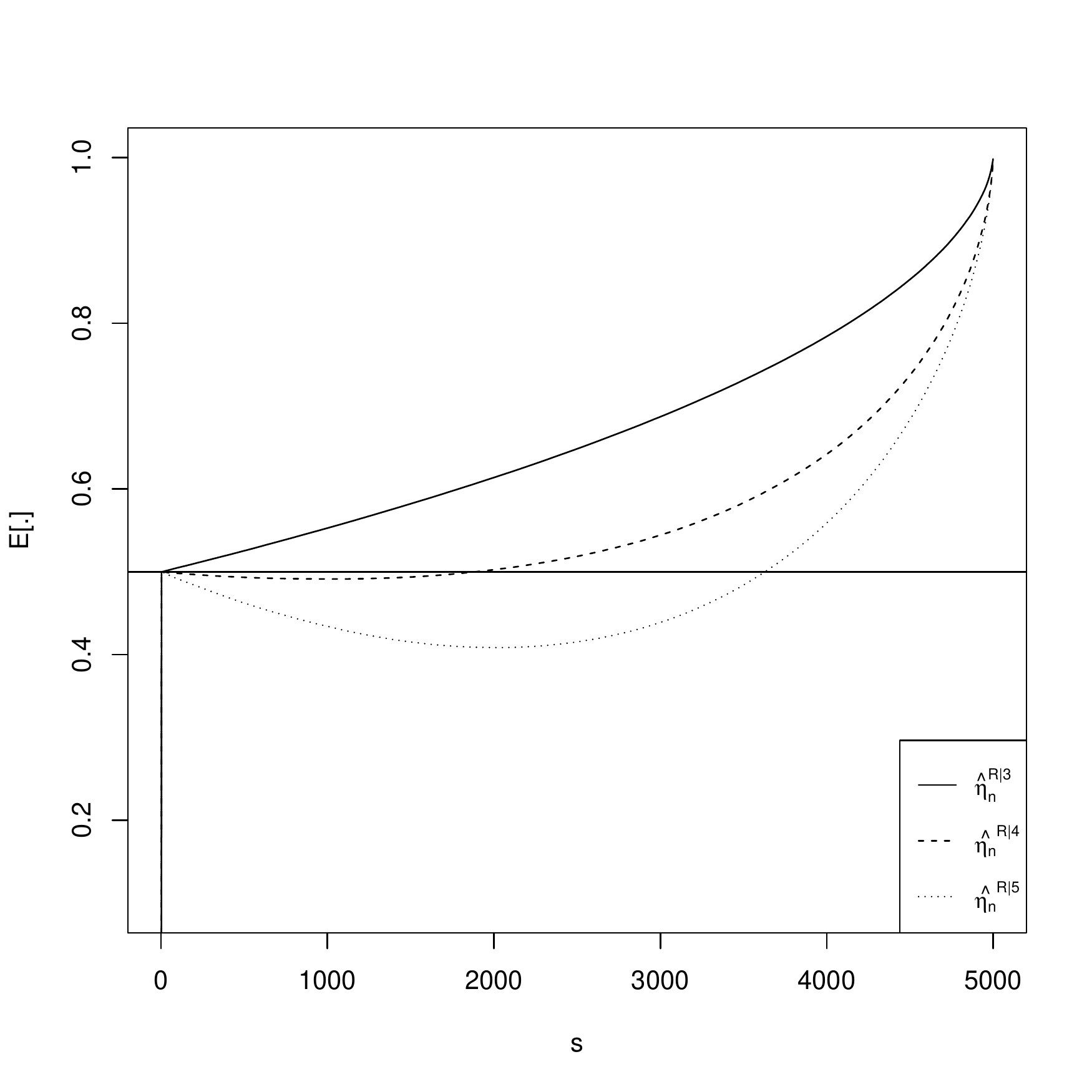}
\includegraphics[scale=0.35]{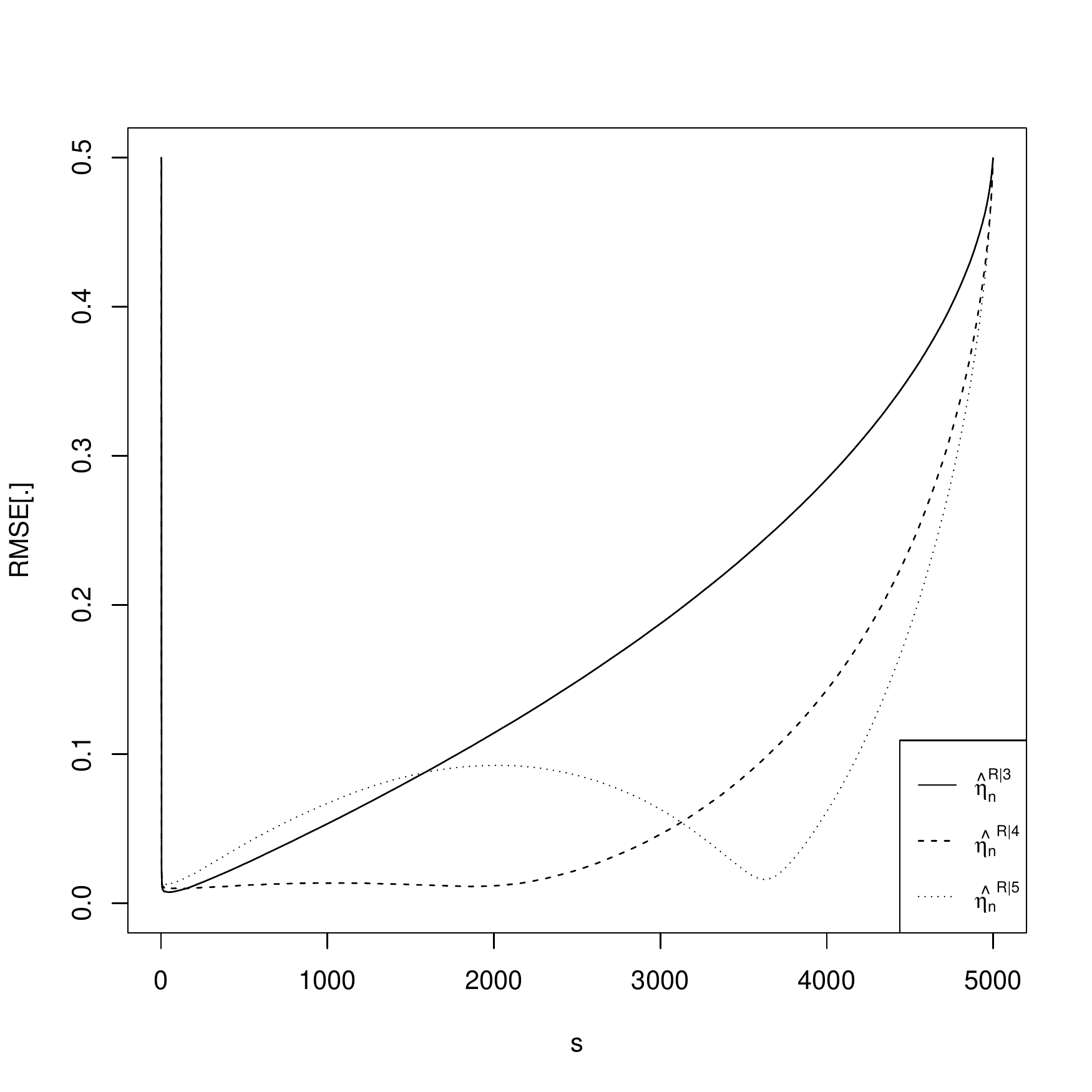}\vspace{-0.5cm}
\caption{Mean values (left) and root mean squared errors (right) of the estimators $\widehat{\eta}^{R|k}_n,$ $k=3,4,5,$ for the MM process in Example 2.1 ($\eta = 0.5$).}
\end{center}
\end{figure}

As we can see, from Figure 4, for the MM process the estimator $\widehat{\eta}^{R|4}_n$ outperforms the estimators $\widehat{\eta}^{R|3}_n$ and $\widehat{\eta}^{R|5}_n,$ improving the results found in Sebastião {\it{et al.}} (2013) \cite{seb1}.  The estimates have a big stability region near the true value $\eta=0.5$ and the RMSE presents a desirable wide ``bathtub'' pattern. For the processes of Examples 2.2 and 2.3 the best results were obtained with $\widehat{\eta}^{R|3}_n.$ The results show that even if we validate condition $\widetilde{D}^{(k)}(u_n)$ for some value $k\geq 3$ we might obtain better estimates for $\eta$ with a runs estimator with a higher value of $k.$ It is therefore convenient to plot sample paths of the runs estimator $\widehat{\eta}^{R|k}_n(X_{n-s:n})$ for different values of $k,$ as a function of $s,$ and analyze the stability regions. Note however that the choice of a too large $k$ may reduce the precision of the new estimators, as can be seen by the values plotted in Figure 4.

\section{Simulation study}

\pg We compare the performance of the several estimators previously presented for a range of different processes, having various dependence structures which give both  $\eta<1$ and $\eta=1$

In this comparison study we include the estimator that naturally arises from relation (\ref{rel_eta_teta}), by re\-pla\-cing $\theta,$ $\tau$ and $\nu$ with corresponding  consistent estimators. Such an estimator has already been considered in Sebastião {\it{et al.}} (2013) \cite{seb1}, where $\theta$ was estimated with  Ferro and Seger's (\cite{ferro}) intervals estimator, $\widehat{\theta}^{FS},$ and $\tau$ and $\nu$ by their empirical counterparts. This estimator is formally defined as follows
\begin{equation}
\widehat{\eta}_n^{EI}=\widehat{\eta}_n^{EI}(u_n):=\widehat{\theta}^{FS}\frac{\displaystyle{\sum_{i=1}^n\indi_{\{X_i\geq u_n\}}}}{\displaystyle{\sum_{i=1}^{n-1}\indi_{\{X_i\leq u_n<X_{i+1}\}}}},\label{estim_rel}
\end{equation}
where
$$\widehat{\theta}^{FS}=\left\{
\begin{array}{ccl}
1\wedge \hat{\theta}^{1}_n(u_n)& , & \max\{T_i:\ 1\leq i\leq N-1\}\leq 2\\
1\wedge \hat{\theta}^{2}_n(u_n) & , &  \max\{T_i:\ 1\leq i\leq N-1\}> 2
\end{array}\right.$$
with $$\hat{\theta}^{1}_n(u_n)=\frac{\displaystyle{2\left(\sum_{i=1}^{N-1}T_i\right)^2}}{\displaystyle{(N-1)\sum_{i=1}^{N-1}T_i^2}}\quad {\textrm{and}}\quad \hat{\theta}^{2}_n(u_n)=\frac{\displaystyle{2\left\{\sum_{i=1}^{N-1}(T_i-1)\right\}^2}}{\displaystyle{(N-1)\sum_{i=1}^{N-1}(T_i-1)(T_i-2)}},$$
denoting $T_i,$ $i=1,\ldots,N-1,$ the observed interexceedance times, {\it{i.e.}}, $T_i=S_{i+1}-S_i,$ where \linebreak $1\leq S_1<\ldots<S_N\leq n$ are the exceedance times.

The different estimation methods are compared via the estimated bias and variance. Furthermore, to analyze the bias-variance trade-off we compare the corresponding root mean squared errors, RMSE. To measure the degree of possible misinterpretation caused by considering the variance/standard error as a criteria for the quality of the estimator we consider the mean error to standard errors $\textrm{MESE}=\textrm{Bias}/\sqrt{\textrm{Var}}$ (see Frahm {\it{et al.}} (2005) \cite{frah} for further details).

The values for the statistical quantities considered were obtained with 5000 independent replications of the estimation procedures for samples of size $n=5000$. Since all estimators of $\eta$ depend on a threshold, we considered all simulated values for each estimator computed at thresholds corresponding to the 90\% and 95\% empirical quantiles.  To reduce sampling variation, for each process all the estimators are evaluated using the same simulated data sets.

For the blocks estimators we compare different block sizes, found to be the best choices from independent simulations to those reported here and from the results obtained in Tables 1-3.

The performance of the estimators for the processes given in Examples 2.1, 2.2 and 2.3 are presented, respectively, in Tables 4, 5 and 6.

\begin{table}[!htb]\begin{center}
 \caption{Simulation results for the estimators of $\eta$ with the MM process of Example 2.1.}\label{c4:tab:optim:mm}\vspace{0.2cm}
 \begin{tabular}{l|cccc|cccc}\hline\hline
  Estimator & Bias &Var &RMSE&MESE&Bias &Var &RMSE&MESE \\ \cline{2-9}
 &\multicolumn{4}{c|}{$q_{0.90}$}&\multicolumn{4}{c}{$q_{0.95}$}\\\hline
 $\widehat{\eta}_n^B \ (r_n=5)$ & $\ \ 0.1964^-$ & 0.00033 & $\ 0.1972^-$ & $10.8028^-$ & $\ 0.1983^-$ & 0.00069 & $\ 0.2000^-$ &
$\ 7.5571^-$\\ [0.3cm]
 $\widehat{\eta}_n^B \ (r_n=10)$ & 0.0561 & 0.00025 & 0.0583 & 3.5137 & 0.0770 & 0.00049 & 0.0801 &
$\ 3.4779^-$\\ [0.3cm]
 $\widehat{\eta}_n^B \ (r_n=15)$ & $-0.0126^+ \ $ & 0.00022 & $\ 0.0196^+$ & $\ 0.8459^+$ & 0.0250 & 0.00044 & 0.0326 &
1.1946\\ [0.3cm]
 $\widehat{\eta}_n^{dj} \ (r_n=50)$ & 0.1195 & 0.00404 & 0.1353 & 1.8797 & 0.0717 & 0.00228 & 0.0862
& 1.5002\\ [0.3cm]
 $\widehat{\eta}_n^{dj} \ (r_n=60)$ & 0.1179 & 0.00602 & 0.1411 & 1.5204 & 0.0667 & $\ 0.00286^-$ & 0.0855
& 1.2471\\ [0.3cm]
 $\widehat{\eta}_n^{dj} \ (r_n=70)$ & 0.1172 & $\ 0.00898^-$ & $\ 0.1507^-$ & 1.2366 & 0.0660 & 0.00361 & 0.0893
& 1.0990\\ [0.3cm]
 $\widehat{\eta}_n^{sl} \ (r_n=50)$ & 0.1156 & 0.00202 & 0.1241 & 2.5719 & 0.0705 & 0.00129 & 0.0791
& 1.9663\\ [0.3cm]
 $\widehat{\eta}_n^{sl} \ (r_n=60)$ & 0.1124 & 0.00283 & 0.1243 & 2.1139 & 0.0660 & 0.00164 & 0.0774
& 1.6305\\ [0.3cm]
 $\widehat{\eta}_n^{sl} \ (r_n=70)$ & 0.1117 & 0.00396 & 0.1282 & 1.7748 & 0.0633 & 0.00205 & 0.0778
& 1.3993\\ [0.3cm]
 $\widehat{\eta}_n^{R|3}$ & 0.0255 & $\ 0.00004^+$ & 0.0262 & $\ 4.0038^-$ & $\ 0.0126^+$ & $\ 0.00004^+$ & $\ 0.0141^+$ & 1.9755\\
[0.3cm]
 $\widehat{\eta}_n^{R|4}$ & $-0.0067^+\ $ & $\ 0.00009^+$ & $\ 0.0117^+$ & $\ 0.6892^+$ & $-0.0039^+\ $ & $\ 0.00009^+$ & $\ 0.0104^+$ & $\ 0.4029^+$\\
[0.3cm]
 $\widehat{\eta}_n^{R|5}$ & $-0.0376\ \ $ & 0.00013 & 0.0394 & 3.2787 & $-0.0202\ \ $ & 0.00014 & 0.0234 &
1.7051\\ [0.3cm]
 $\widehat{\eta}_n^{FS}$ & 0.0506 & 0.00183 & 0.0663 & 1.1827 & 0.0300 & $\ 0.00308^-$ & 0.0631 & $\ 0.5415^+$\\
[0.3cm]
 \hline\hline
 \end{tabular}\\
Best values are ticked with a plus and worst are ticked with a minus.
\end{center}
\end{table}

\begin{table}[!htb]\begin{center}
 \caption{Simulation results for the estimators of $\eta$ with the AR(1) process with $\beta=2$ of Example 2.2.}\vspace{0.2cm}
 \begin{tabular}{l|cccc|cccc}\hline\hline
  Estimator & Bias &Var &RMSE&MESE&Bias &Var &RMSE&MESE \\ \cline{2-9}
 &\multicolumn{4}{c|}{$q_{0.90}$}&\multicolumn{4}{c}{$q_{0.95}$}\\\hline
 $\widehat{\eta}_n^B \ (r_n=5)$ & 0.0881 & 0.00022 & 0.0894 & 5.9123 & 0.1008 & $\ 0.00045^+$ & 0.1030 &
$\ 4.7764^-$\\ [0.3cm]
 $\widehat{\eta}_n^B \ (r_n=10)$ & $-0.0874\ \ $ & 0.00024 & 0.0888 & 5.6748 & 0.0009  & $\ 0.00055^+$ & $\ 0.0235^+$ &
0.0399\\ [0.3cm]
 $\widehat{\eta}_n^B \ (r_n=15)$ & $-0.2139\ \ $ & $\ 0.00015^+$ & 0.2143 & $17.6109^-$ & $-0.0776\ \ $ & $\ 0.00055^+$ & 0.0810
& $\ 3.3176^-$\\ [0.3cm]
 $\widehat{\eta}_n^{dj} \ (r_n=50)$ & 0.1404 & 0.00393 & 0.1538 & 2.2390 & 0.1974 & 0.01331 & 0.2286
& 1.7104\\ [0.3cm]
 $\widehat{\eta}_n^{dj} \ (r_n=60)$ & $-0.0221\ \ $ & 0.00095 & 0.0380 & 0.7167 & $\ 0.2093^-$ & $\ 0.02289^-$ &
$\ 0.2582^-$ & 1.3834\\ [0.3cm]
 $\widehat{\eta}_n^{dj} \ (r_n=70)$ & $-0.1428\ \ $ & $\ 0.00020^+$  & 0.1434 & 10.1522 & 0.2048 & $\ 0.02340^-$ & $\ 0.2556^-$
& 1.3385\\ [0.3cm]
 $\widehat{\eta}_n^{sl} \ (r_n=50)$ & $\ 0.3980^-$ & $\ 0.04007^-$ & $\ 0.4455^-$ & 1.9883 & 0.1890 & 0.00698 & 0.2066
& 2.2622\\ [0.3cm]
 $\widehat{\eta}_n^{sl} \ (r_n=60)$ & $\ 0.3100^-$ & $\ 0.02936^-$ & $\ 0.3542^-$ & 1.8094 & 0.1949 & 0.01014 & 0.2194
& 1.9349\\ [0.3cm]
 $\widehat{\eta}_n^{sl} \ (r_n=70)$ & 0.1953 & 0.01983 & 0.2408 & 1.3870 & $\ 0.2050^-$ & 0.01545 & 0.2397
& 1.6490\\ [0.3cm]
 $\widehat{\eta}_n^{R|3}$ & $\ 0.0006^+$  & 0.00039 & $\ 0.0198^+$ & $\ 0.0279^+$ & 0.0007 & 0.00077 & $\ 0.0278^+$ & 0.0255\\
[0.3cm]
 $\widehat{\eta}_n^{R|4}$ & $\ 0.0004^+$  & 0.00039 & $\ 0.0198^+$ & $\ 0.0193^+$ & $\ 0.0006^+$ & 0.00077 & $\ 0.0278^+$ & $\ 0.0208^+$\\
[0.3cm]
 $\widehat{\eta}_n^{R|5}$ & $-0.0622\ \ $ & 0.00044 & 0.0657 & 2.9574 & $\ 0.0004^+$  & 0.00077 & $\ 0.0278^+$ &
$\ 0.0159^+$\\ [0.3cm]
 $\widehat{\eta}_n^{FS}$ & 0.2424 & 0.00040 & 0.2432 & $\ 12.0808^-$ & 0.1566 & 0.00394 & 0.1687 & 2.4965\\
[0.3cm]
 \hline\hline
 \end{tabular}\\
Best values are ticked with a plus and worst are ticked with a minus.
\end{center}
\end{table}

\begin{table}[!htb]\begin{center}
 \caption{Simulation results for the estimators of $\eta$ with the MAR(1) process with $\alpha=0.9$ of Example 2.3.}\vspace{0.2cm}
 \begin{tabular}{l|cccc|cccc}\hline\hline
 Estimator & Bias &Var &RMSE&MESE&Bias &Var &RMSE&MESE \\ \cline{2-9}
 &\multicolumn{4}{c|}{$q_{0.90}$}&\multicolumn{4}{c}{$q_{0.95}$}\\\hline
 $\widehat{\eta}_n^B \ (r_n=5)$ & $-0.0019^+\ $ & $\ 0.00004^+$ & $\ 0.0065^+$ & $\ 0.3006^+$ & $-0.0008^+ \ $ & $\ 0.00003^+$ & $\ 0.0056^+$ &
$\ 0.1420^+ $\\ [0.3cm]
 $\widehat{\eta}_n^B \ (r_n=10)$ & $-0.0100\ \ $ & 0.00020 & 0.0174 & 0.7001 & $-0.0050\ \ $ & 0.00019 & 0.0148 &
0.3583\\ [0.3cm]
 $\widehat{\eta}_n^B \ (r_n=15)$ & $-0.0232\ \ $ & 0.00046 & 0.0317 & 1.0804 & $-0.0112\ \ $ & 0.00043 & 0.0236
& 0.5397\\ [0.3cm]
 $\widehat{\eta}_n^{dj} \ (r_n=50)$ & $\ 0.3304^-$ & 0.01156 & $\ 0.3474^-$ & 3.0723 & 0.2586 & 0.01341 & 0.2834
& 2.2335\\ [0.3cm]
 $\widehat{\eta}_n^{dj} \ (r_n=60)$ & 0.2986 & 0.01363 & 0.3206 & 2.5580 & 0.2248 & 0.01387 & $\ 0.2538^-$
& 1.9086\\ [0.3cm]
 $\widehat{\eta}_n^{dj} \ (r_n=70)$ & 0.2759 & $\ 0.01526^-$ & 0.3023 & 2.2332 & 0.2066 & $\ 0.01429^-$ & 0.2386
& 1.7283\\ [0.3cm]
 $\widehat{\eta}_n^{sl} \ (r_n=50)$ & $\ 0.3270^-$ & 0.00591 & $\ 0.3359^-$ & $\ 4.2524^-$ & $\ 0.2552^-$ & 0.00594 & $\ 0.2666^-$
& $\ 3.3119^-$\\ [0.3cm]
 $\widehat{\eta}_n^{sl} \ (r_n=60)$ & 0.2949 & 0.00687 & 0.3063 & $\ 3.5595^-$ & 0.2241 & 0.00626 & 0.2377
& 2.8332\\ [0.3cm]
 $\widehat{\eta}_n^{sl} \ (r_n=70)$ & 0.2726 & 0.00806 & 0.2870 & 3.0372 & 0.2028 & 0.00695 & 0.2192
& 2.4321\\ [0.3cm]
 $\widehat{\eta}_n^{R|3}$ & $-0.0013^+ \ $ & $\ 0.00003^+$ & $\ 0.0055^+$ & $\ 0.2518^+$ & $-0.0009^+ \ $ & $\ 0.00003^+$ & $\ 0.0060^+$ & $\ 0.1451^+$\\
[0.3cm]
 $\widehat{\eta}_n^{R|4}$ & $-0.0035\ \ $ & 0.00007 & 0.0091 & 0.4140 & $-0.0021\ \ $ & $\ 0.00009^+$ & 0.0095 & 0.2250\\
[0.3cm]
 $\widehat{\eta}_n^{R|5}$ & $-0.0066\ \ $ & 0.00014 & 0.0135 & 0.5666 & $-0.0036\ \ $ & 0.00014 & 0.0125 &
0.3019\\ [0.3cm]
 $\widehat{\eta}_n^{FS}$ & 0.0324 & $\ 0.01645^-$ & 0.1323 & 0.2523 & 0.0426 & $\ 0.02779^-$ & 0.1720 & 0.2558\\
[0.3cm]
 \hline\hline
 \end{tabular}\\
Best values are ticked with a plus and worst are ticked with a minus.
\end{center}
\end{table}

For a better comparison of the blocks estimators $\widehat{\eta}^{B}_n$ with $r_n=5,10,15$ and the runs estimators $\widehat{\eta}^{R|k}_n,$ $k=3,4,$  which had the best performances in the previous simulations,  we plot in Figure 5 the estimated mean values and RMSE for levels given by the $(s+1)$th, $s=1,...,n-1,$ top order statistics.

\begin{figure}[!ht]
\begin{center}
\includegraphics[scale=0.3]{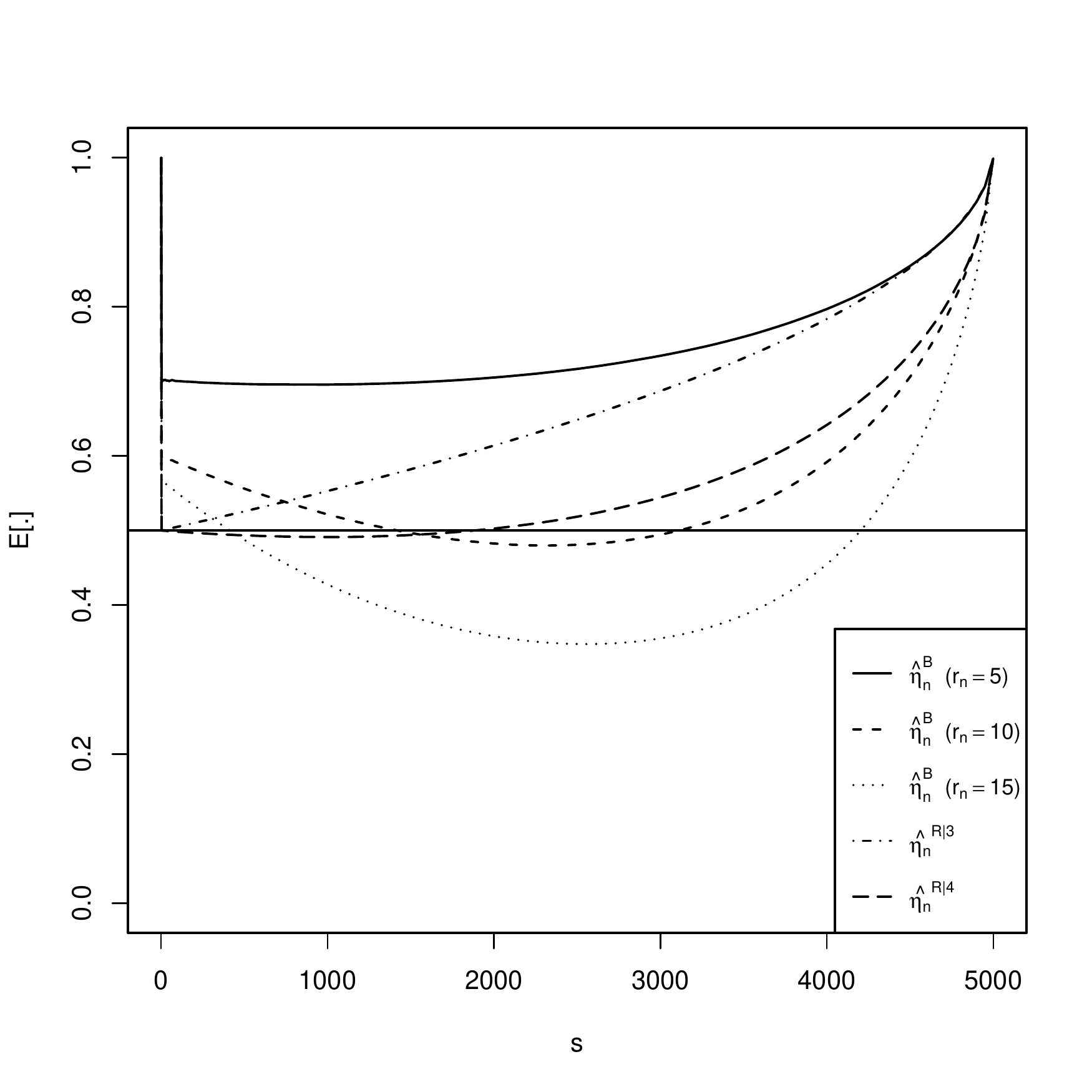} \includegraphics[scale=0.3]{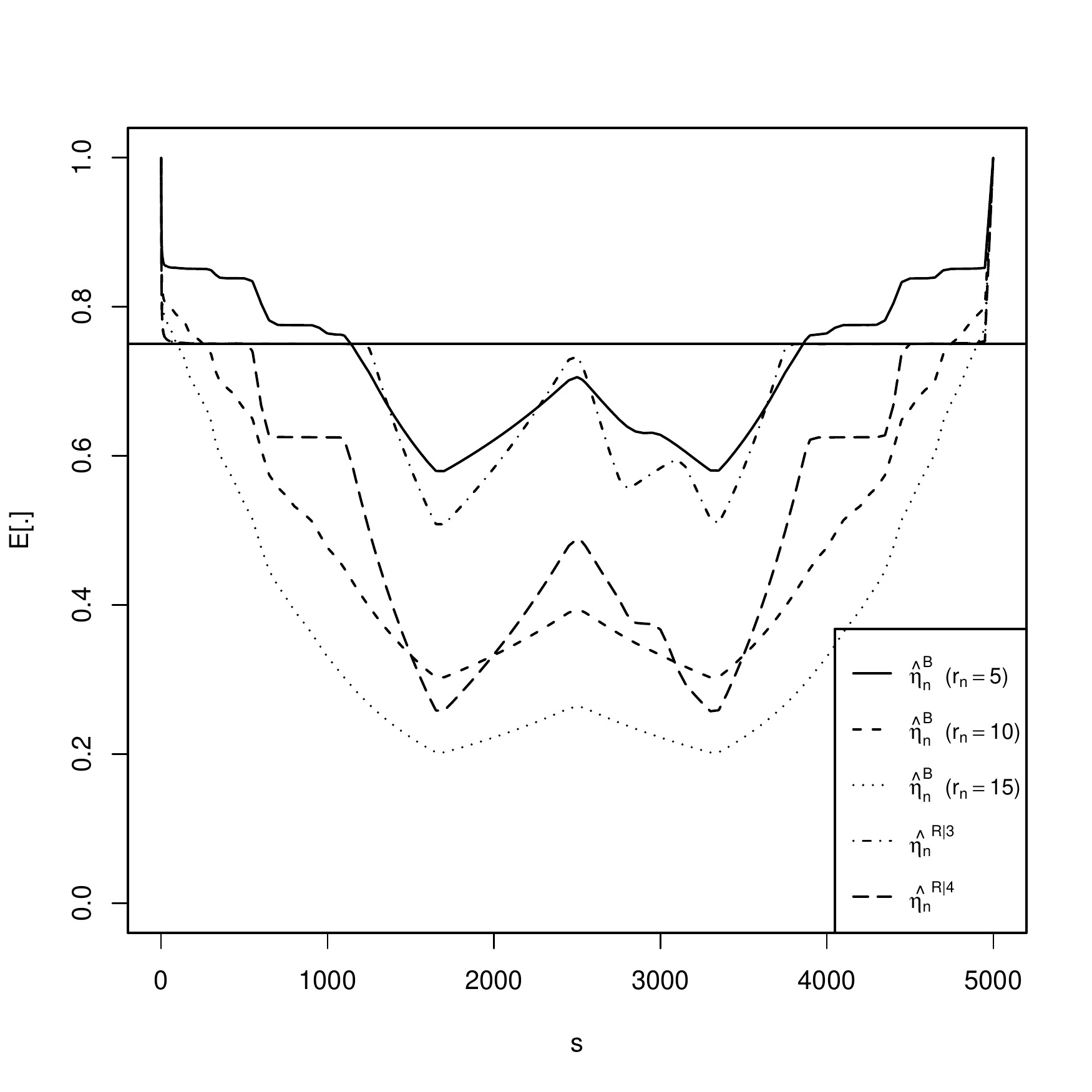} \includegraphics[scale=0.3]{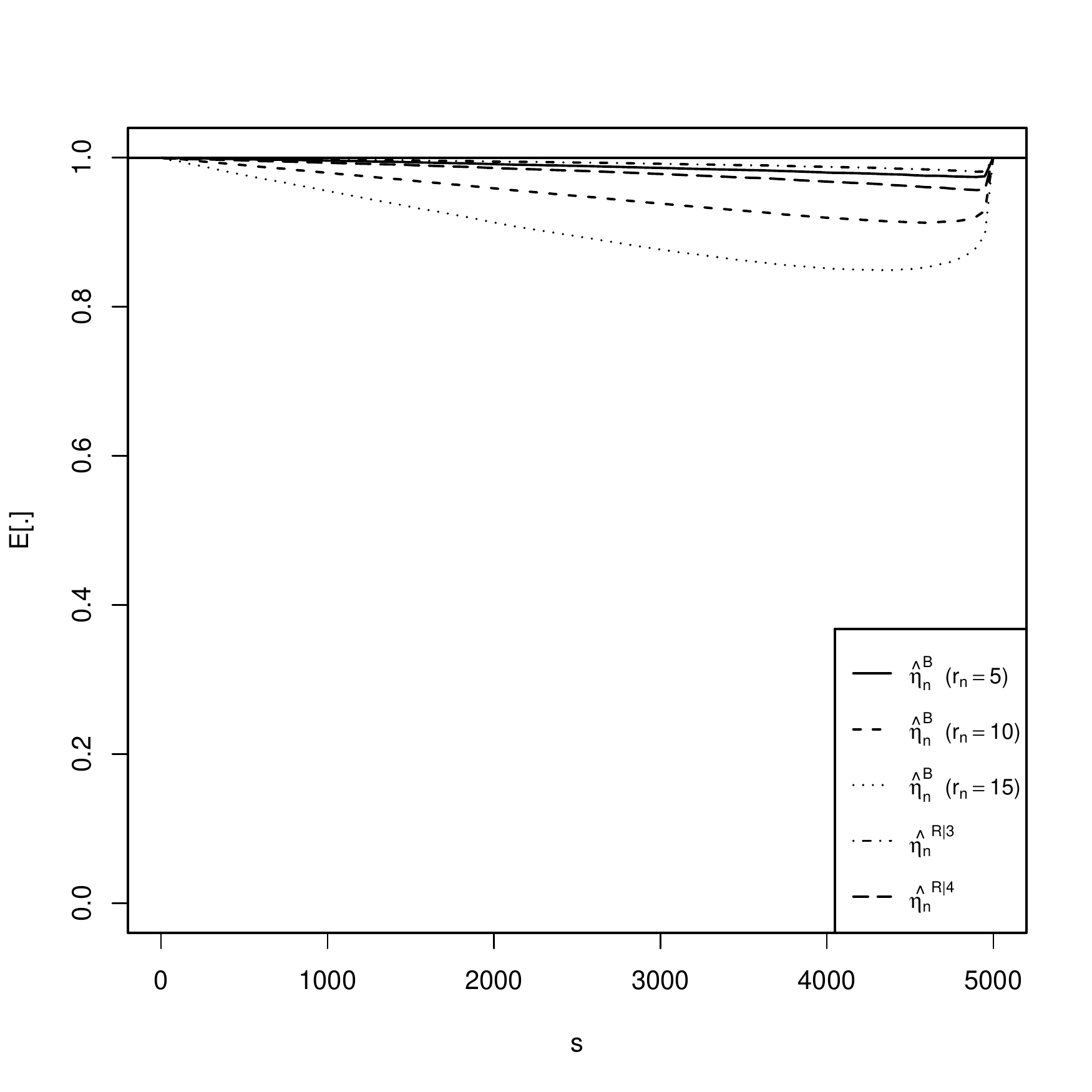}\\
\includegraphics[scale=0.3]{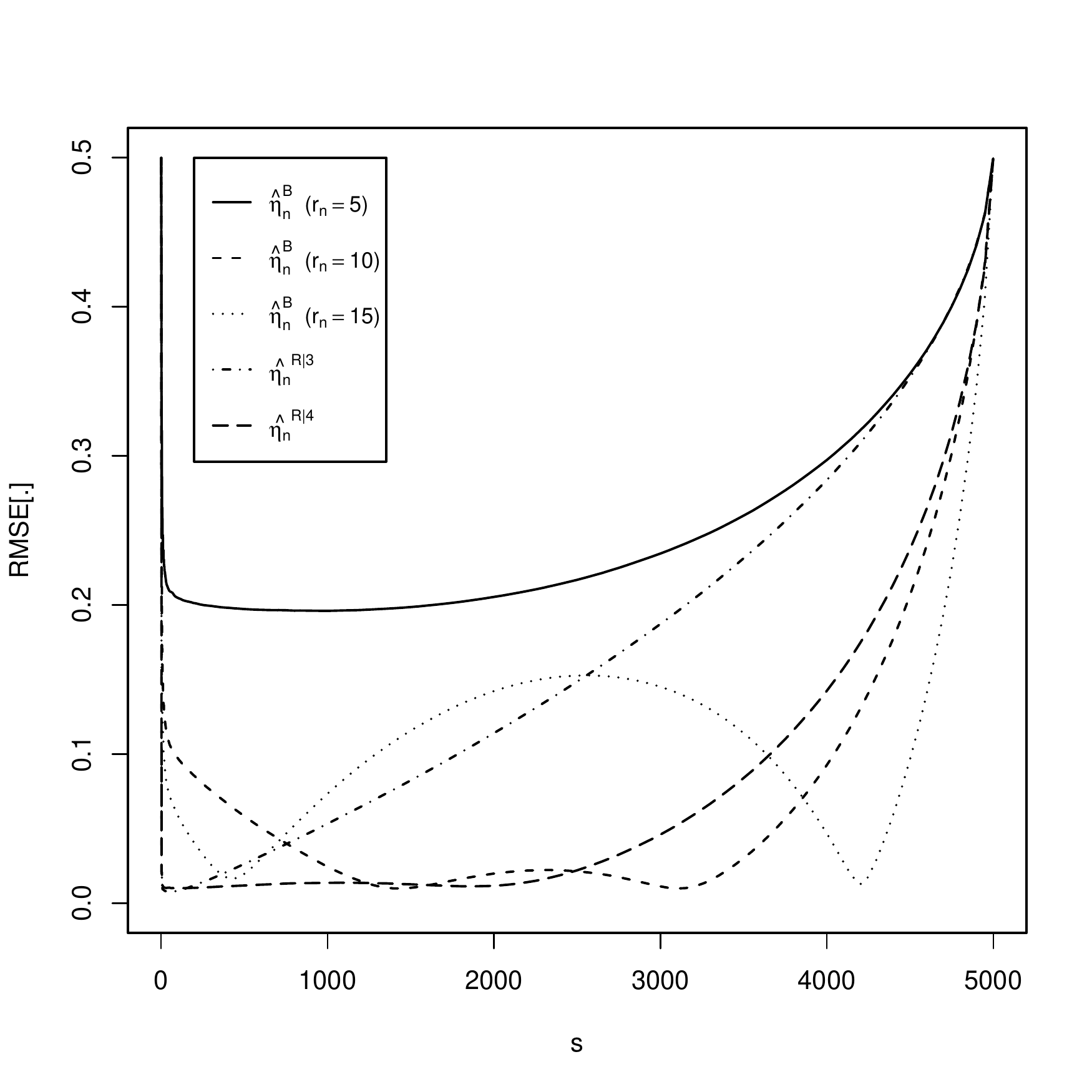} \includegraphics[scale=0.3]{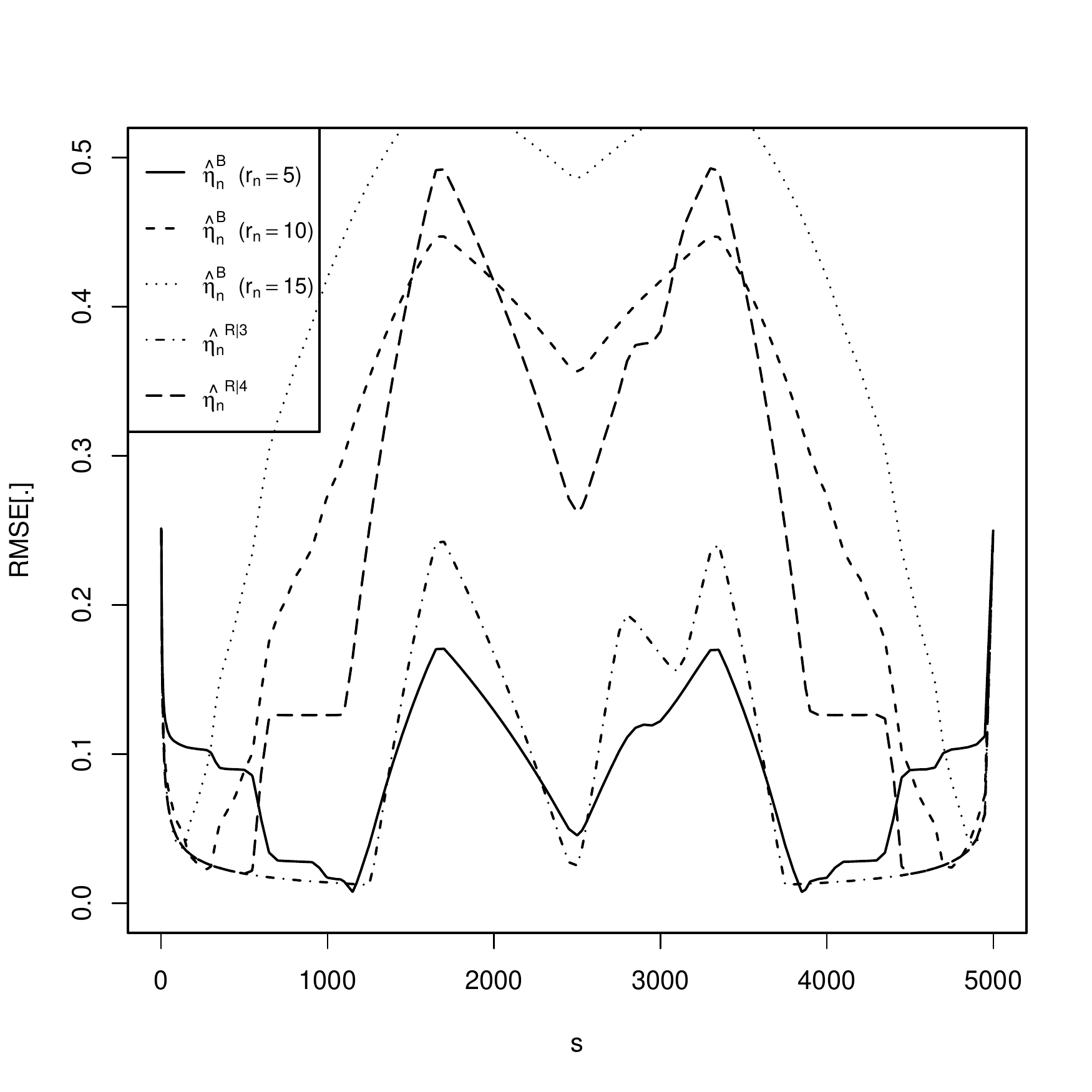}\includegraphics[scale=0.3]{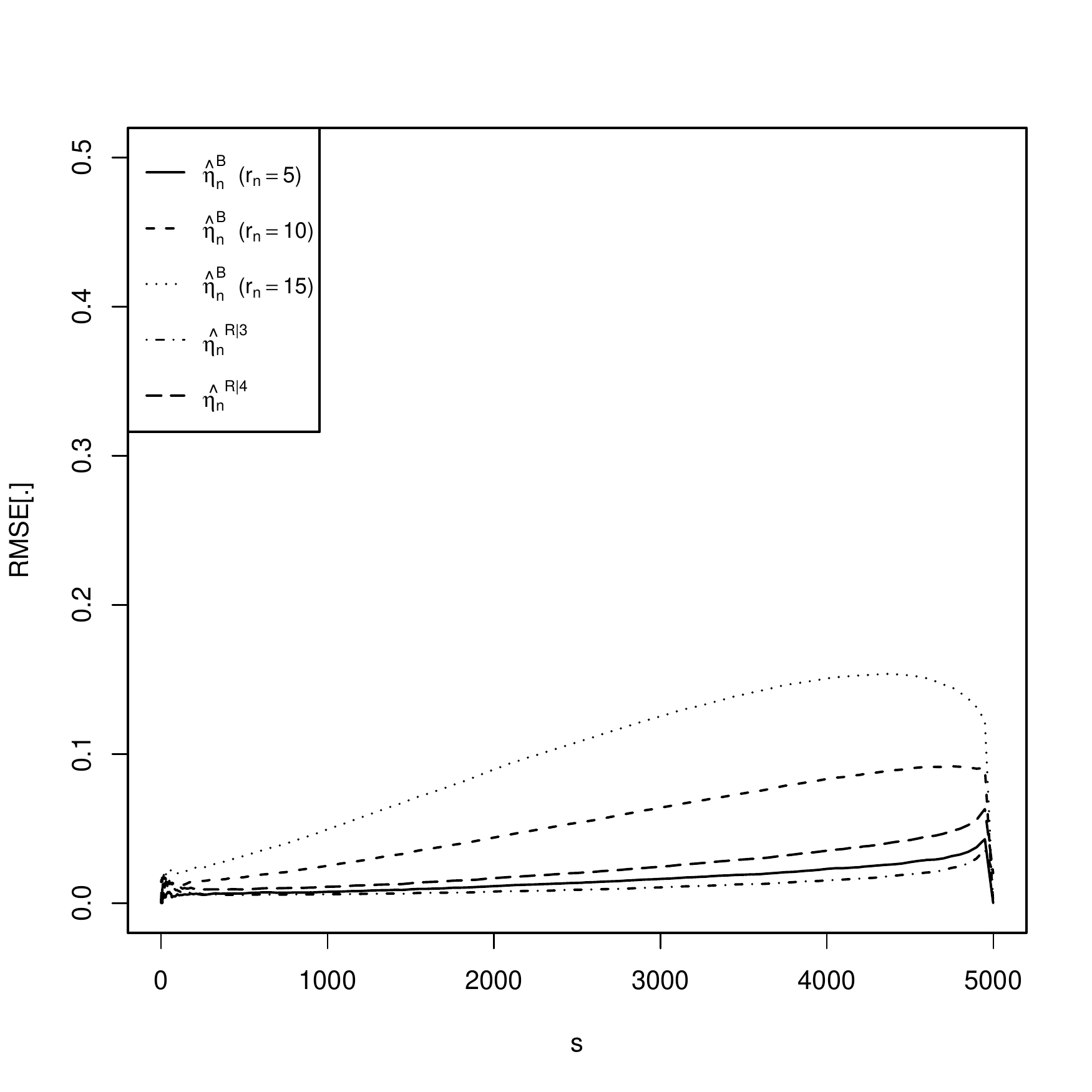}\vspace{-0.5cm}
\caption{From left to right, mean values (top) and root mean squared errors (bottom) of the estimators $\widehat{\eta}^{B}_n$ with $r_n=5,10,15$ and  $\widehat{\eta}^{R|k}_n,$ $k=3,4,$ for the MM process in Example 2.1 ($\eta = 0.5$), the AR(1) process with $\beta=2$ in Example 2.2 ($\eta=0.75$) and the MAR(1) process with $\alpha=0.9$ in Example 2.3 ($\eta=1$). }
\end{center}
\end{figure}


\vspace{0.3cm}\noindent {\bf{Some overall conclusions:}}
\begin{itemize}
\item In general, all estimation methods have small sample bias. For the three processes the blocks estimators $\widehat{\eta}_n^B,$ $r_n=5,$ 10,15 and the runs estimators $\widehat{\eta}_n^{R|k},$ $k=3,4,5,$ lead to the smallest absolute sample biases.  The biggest sample biases are obtained with the sliding blocks estimators $\widehat{\eta}_n^{sl}$ and the disjoint blocks estimators $\widehat{\eta}_n^{dj}.$ The latter has a bit worse performance, except for the AR(1) process at the 90\% threshold where the sliding blocks estimator and the naive estimator $\widehat{\eta}_n^{FS}$ perform worst. 
\item The sample variances for all estimators and processes are smaller than 0.05. In the group of the smallest observed variance values we find once again the the blocks estimators $\widehat{\eta}_n^B,$ $r_n=5,$ 10,15 and the runs estimators $\widehat{\eta}_n^{R|k},$ $k=3,4,5.$ The conclusions drawn in Remark 7 can be pictured in Tables 4, 5 and 6. When $\eta$ is on the boundary of the parameter space there is no clear overall change in efficiency of the estimators with threshold level.
\item The conclusions previously drawn for the bias and variance carry over to the RMSE, for obvious reasons.  For the MAR(1) process, with $\eta=1,$ we obtain the smallest values of the RMSE for all estimators, in fact the blocks estimator $\widehat{\eta}_n^B,$ $r_n=5,$ and the runs estimator $\widehat{\eta}_n^{R|3}$ have a RMSE of 0.006.
\item A large sample bias relative to the sample variance translates into a large MESE, which is true for the blocks estimators  $\widehat{\eta}_n^B,$ $r_n=5,10,$ with the MM and the AR(1) processes. The runs estimators present in the majority of the cases values smaller than 1, which indicates that the true upcrossings index $\eta$ lies within $1\sigma$ confidence band, where $\sigma$ denotes the standard error. 
\item All but the blocks estimator  $\widehat{\eta}_n^B$ have a better performance the higher the threshold used, since in essence they are estimating the non-upcrossing or non-exceedance of the threshold.
\item Estimators obtained from asymptotic characterizations of the upcrossings index have a better performance than estimators obtained from the relation with the extremal index.
\item The good performance of the blocks estimator $\widehat{\eta}_n^B$ makes it the best alternative for the runs estimators  $\widehat{\eta}_n^{R|k}$ which need the validation of local dependence conditions. This statement is reinforced by Figure 5, where we see that for an adequate block size the estimated mean of the blocks estimator presents a well defined stability region near the true value and the RMSE has a wide ``bathtub'' shape.  
\end{itemize}



\section{Application to financial data}

\pg We consider the performance of the estimators under study which have previously better performed, {\it{i.e.}} the blocks estimator $\widehat{\eta}_n^B$ and the runs estimators $\widehat{\eta}_n^{R|k},$ $k\geq 3,$ when applied to the analysis of the log-returns associated with the daily closing prices of the German stock market index DAX, collected from  3 January 2000 to 29 June 2012. This series was analysed previously by Sebastião {\it{et al.}} (2013) \cite{seb1} and is plotted in Figure 6: DAX daily closing prices over the mentioned period, $x_t$, and the log-returns, $100 \times (\ln x_t - \ln x_{t-1})$, the data to be analyzed, after removing null log-returns ($n=3177$).

\begin{figure}[!ht]
\begin{center}
\includegraphics[scale=0.45]{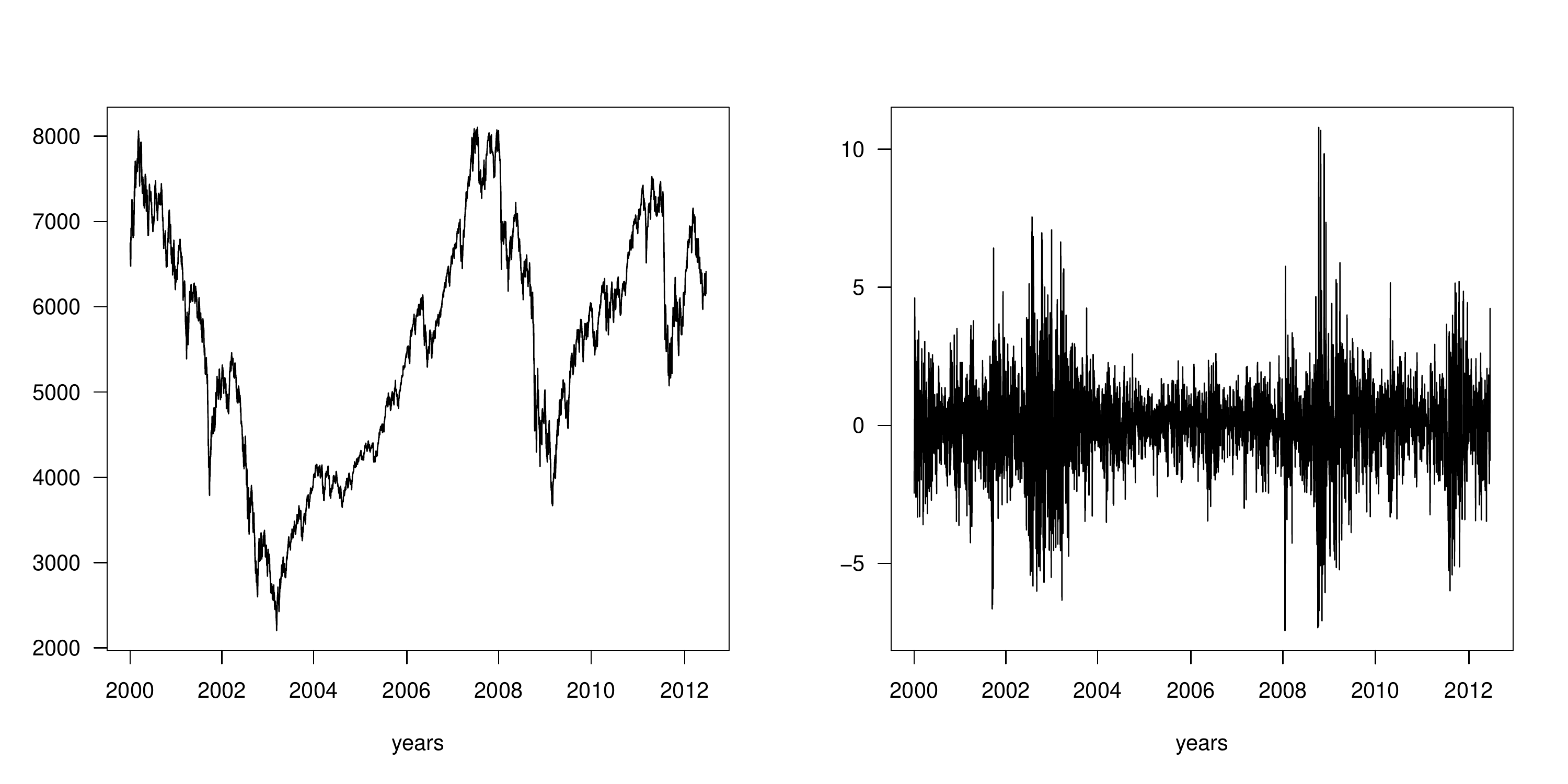}\vspace{-0.5cm}
\caption{DAX daily closing prices (left) and daily log-returns  (right), of DAX from 2000 to 2012,  with 3177 observations (successive equal prices excluded).}
\end{center}
\end{figure}

As stated in Sebastião {\it{et al.}} (2013) \cite{seb1},  Klar et al. (2012) \cite{klar} have analyzed the DAX German stock market index time series and concluded that the GARCH(1,1) process with Student-$t$ distributed innovations is a good model to describe these data. The simulated values of the anti-$\widetilde{D}^{(k)}(u_n),$  $k=3,4,5,6$ proportions found in Figure 3 for the GARCH(1,1) process and the conclusions there taken led us to consider the runs estimators $\widehat{\eta}^{R|k}$ with $k=4,$ 5 and $6.$ For the blocks estimator, the choice of $r_n=5,$ 7 and 10 seemed adequate for the sample size considered.

Clusters of upcrossings of high returns can bring serious risk to financial investors so estimates of the upcrossings index for the log-returns are plotted against $s\geq 1$ in Figure 7, in a linear scale and in logarithmic scale to better observe estimates at high levels. Note that Figure 7 also enhances the effect that the block size has on the estimates.

\begin{figure}[!ht]
\begin{center}
\includegraphics[scale=0.35]{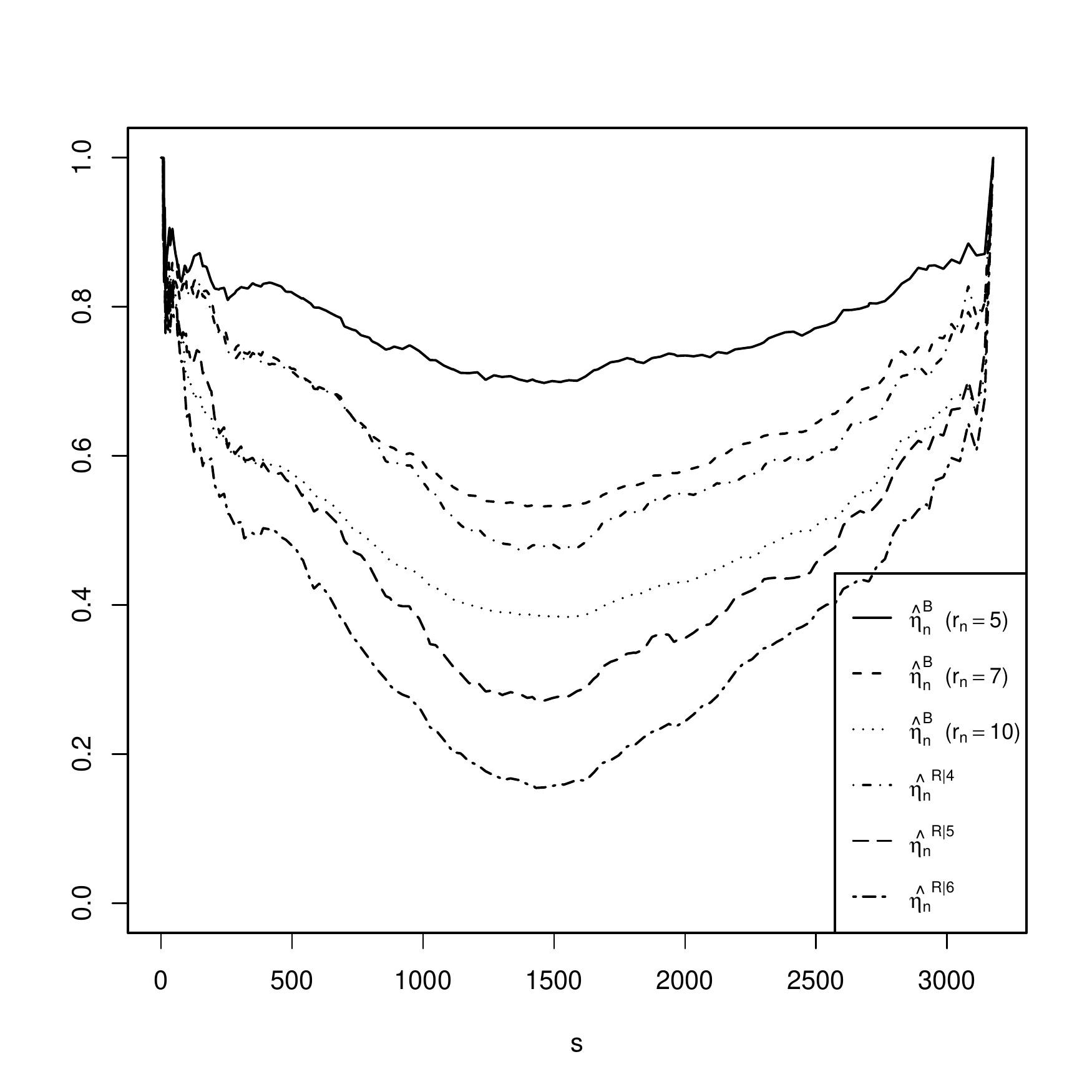} \includegraphics[scale=0.35]{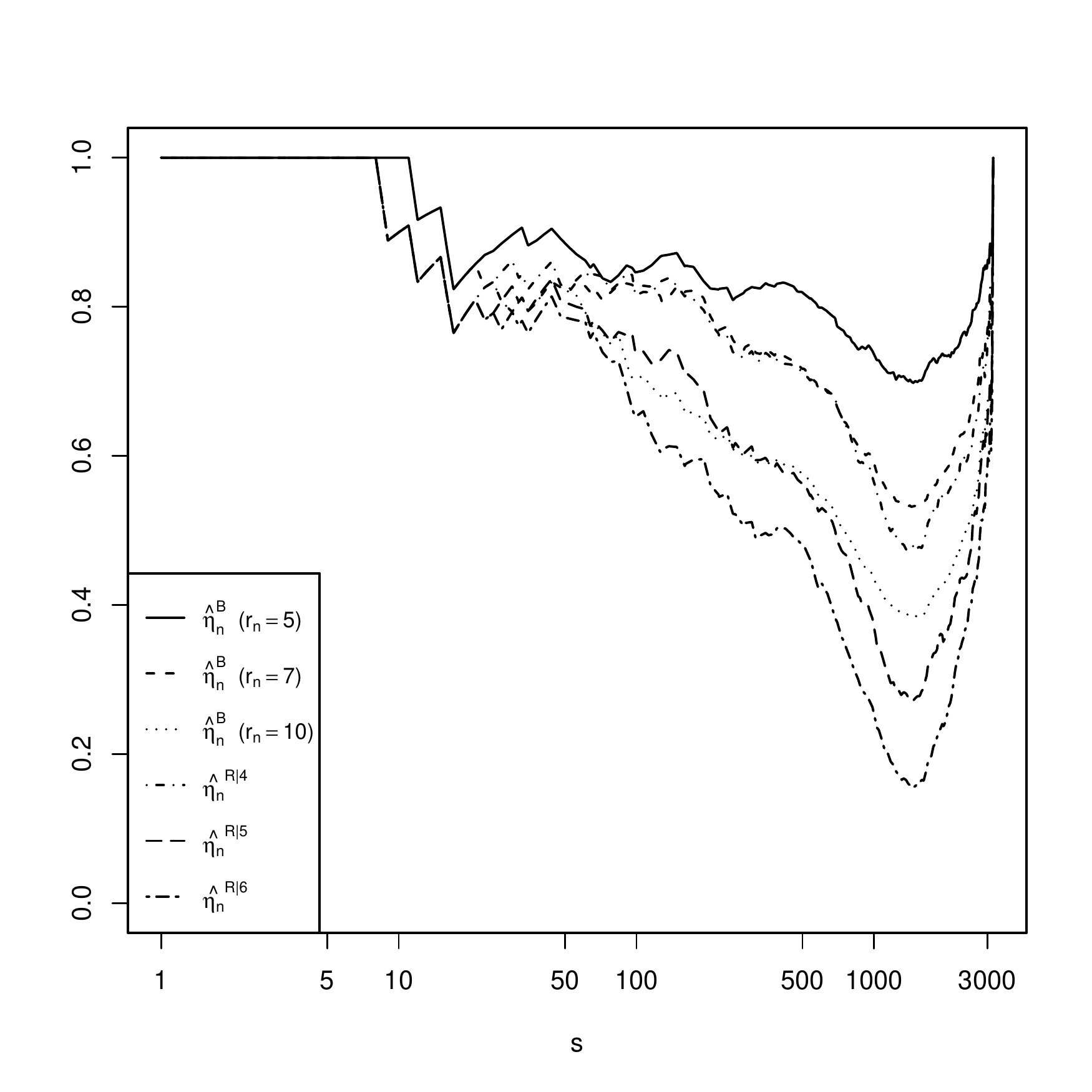}\vspace{-0.5cm}
\caption{Estimates of the upcrossings index for the daily DAX log returns plotted against $s\geq1$, in a linear scale (left) and in a logarithmic scale (right).}
\end{center}
\end{figure}

The sample paths of all estimators have a stability region around the value $\eta=0.8,$ which can be identified in the plot of the estimates in a logarithmic scale of Figure 7. The longest stability region is obtained with the blocks estimator for $r_n=7.$ The higher value of $\eta$ obtained with the runs estimator $\widehat{\eta}^{R|4}$ comparatively to the ones obtained with the runs estimators $\widehat{\eta}^{R|5}$ and $\widehat{\eta}^{R|6}$ can be justified by the fact that condition $\widetilde{D}^{(4)}$ does not hold and therefore runs of upcrossings are over evaluated.

A value of $\eta$ significatively less than 1 shows certain short range dependence reflected in the clustering of high level upcrossings. This may be interpreted as a pattern in the occurrence of the extreme values, that is, once a large profit in the asset return has occurred it will be followed by a loss (negative return) and we can expect a period of large profits followed by losses (values crossing some threshold). The average length of this period being the inverse of the upcrossings index.

\section{Conclusions}

\pg New blocks estimators for the upcrossings index $\eta$, a measure of the clustering of upcrossings of high thresholds, have been proposed in this paper. These estimators have the advantage of not requiring the validation of a local dependence condition, contrarily to the known runs estimator of the upcrossings index. The study of the properties of these blocks estimators, namely their  consistency and asymptotic normality, has been carried out in Section 2. The family of runs estimators of $\eta$ has been enlarged and an empirical way of checking local dependence conditions that control the clustering of upcrossings  provided.

The estimators proposed in this paper improve substantially on existing methods of estimating the upcrossings index.  Nevertheless, the simulations performed, in Sections 3 and 4, tell us that there is still space for improvements. In future,  a discussion on the optimal choice of the key parameter $r_n$ has to be considered as well as bias reduction techniques, namely the generalized Jackknife methodology. \vspace{0.5cm}

\begin{appendix}\setcounter{equation}{0}

\section*{Appendix A: Proofs for Section 2}\setcounter{section}{0}
\refstepcounter{section}

\subsection{Proof of Theorem \ref{teo:blocos1}}

Let us consider the sample $X_1,\ldots, X_{[n/c_n]},$ divided into $k_n/c_n$ disjoint blocks. We can then apply the arguments used in Lemma 2.1 of Ferreira (2006) \cite{fer1} and conclude that
\begin{equation}
P(\widetilde{N}_{[n/c_n]}(v_n)=0)\sim P^{k_n/c_n}(\widetilde{N}_{r_n}(v_n)=0).\label{eq:bloco2}
\end{equation}
Now, from the definition of the upcrossings index $\eta$ we have  $P(\widetilde{N}_{[n/c_n]}(v_n)=0)\xrightarrow [n\rain]{} e^{-\eta \nu},$ hence applying logarithms on both sides of  (\ref{eq:bloco2}), (\ref{eq:bloco1}) follows immediately.\vspace{0.3cm}

For the conditional mean number of upcrossings in each block, we have from (\ref{eq:bloco1}) and the definition of the thresholds $v_n$
\begin{eqnarray*}
E[\widetilde{N}_{r_n}(v_n)\ |\ \widetilde{N}_{r_n}(v_n)>0]&=&\frac{\sum_{j\geq 1} j\times P(\widetilde{N}_{r_n}(v_n)=j,\ \widetilde{N}_{r_n}(v_n)>0)}{P(\widetilde{N}_{r_n}(v_n)>0)}\\[0.3cm]
&=&\frac{E[\widetilde{N}_{r_n}(v_n)]}{P(\widetilde{N}_{r_n}(v_n)>0)}=\frac{r_nP(X_1\leq v_n<X_2)}{P(\widetilde{N}_{r_n}(v_n)>0)}\\[0.3cm]
&\sim& \frac{r_nc_n\nu/n}{c_n\eta \nu/k_n} ,
\end{eqnarray*} which completes the proof since $k_nr_n\sim n.$\hfill $\square$

\subsection{Proof of Theorem \ref{teo:consit:etaB:2}}

It suffices to show that
\begin{equation}
\lim_{n\to +\infty} \sum_{j=1}^{+\infty} j\widetilde{\pi}_n(j;v_n)=\eta^{-1}\hspace{2cm}\label{p1}
\end{equation}
and
\begin{equation}
\sum_{j=1}^{+\infty}j(\widehat{\widetilde{\pi}}_n(j;v_n)-\widetilde{\pi}_n(j;v_n))\xrightarrow [n\rain]{P}0.
\label{p2}
\end{equation}

(\ref{p1}) follows immediately from Theorem 2.1 since
$$\lim_{n\to +\infty} \sum_{j=1}^{+\infty} j\widetilde{\pi}_n(j;v_n)=\lim_{n\to +\infty}E[\widetilde{N}_{r_n}(v_n)\ |\ \widetilde{N}_{r_n}(v_n)>0]=\eta^{-1}.$$

To show (\ref{p2}), lets start by noting that Theorem 2.2 of Hsing (1991) \cite{hsing1} holds for $\widehat{\widetilde{\pi}}_n(j;v_n)$ and $\widetilde{\pi}_n(j;v_n),$ if $Y_{n1},$ the number of exceedances of $v_n$  in a block of size $r_n,$ is replaced by $\widetilde{N}_{r_n}(v_n).$ The proof now follows from considering $\alpha=1$ and $T(j)=j\indi_{\{j\geq 1\}}$ in this theorem for $\widetilde{N}_{r_n}(v_n)$ and verifying that conditions (a), (c) and (d) of this theorem follow from the assumptions of Theorem 2.1, (2) and (3), (b) follows from Theorem 2.1 and (e) from the fact that $$\lim_{n\to +\infty}\frac{k_n}{c_n} E[\widetilde{N}_{r_n}(v_n)]=\lim_{n\to +\infty}\frac{k_n}{c_n} r_nP(X_\leq v_n <X_2)=\nu.$$ \hfill $\square$

\subsection{Proof of Theorem \ref{prop:normZ}}

From Cramer-Wald's Theorem we know that a necessary and sufficient condition for (\ref{conv_normal}) to hold is that for any $a, b\in \R$
\begin{eqnarray}
\lefteqn{\hspace{-3cm}c_n^{-1/2}\left(a\displaystyle{\sum_{i=1}^{k_n}\left(\widetilde{N}_{r_n}^{(i)}(v_n)-
E[\widetilde{N}_{r_n}^{(i)}(v_n)]\right)}+b \displaystyle{\sum_{i=1}^{k_n}\left(\indi_{\{\widetilde{N}_{r_n}^{(i)}(v_n)>0\}}-
E[\indi_{\{\widetilde{N}_{r_n}^{(i)}(v_n)>0\}}]\right)}\right)}\nonumber\\
&\xrightarrow [n\rain]{d}& \mathcal{N}
(0,\ \nu(a^2\eta \sigma^2+2ab+b^2\eta)).\label{conv_normal2}
\end{eqnarray}
We shall therefore prove (\ref{conv_normal2}). For this, lets consider $$U_{r_n-l_n}^{(i)}(v_n)=\sum_{j=(i-1)r_n+1}^{ir_n-l_n}\indi_{\{X_j\leq v_n<X_{j+1}\}},\qquad V_{l_n}^{(i)}(v_n)=\widetilde{N}_{r_n}^{(i)}(v_n)-U_{r_n-l_n}^{(i)}(v_n),\quad 1\leq i\leq k_n$$ and note that
\begin{eqnarray}
\lefteqn{a\sum_{i=1}^{k_n}\left(\widetilde{N}_{r_n}^{(i)}(v_n)-
E[\widetilde{N}_{r_n}^{(i)}(v_n)]\right)+b \sum_{i=1}^{k_n}\left(\indi_{\{\widetilde{N}_{r_n}^{(i)}(v_n)>0\}}-
E[\indi_{\{\widetilde{N}_{r_n}^{(i)}(v_n)>0\}}]\right)}\nonumber\\[0.3cm]
&=&a\sum_{i=1}^{k_n}\left(U_{r_n-l_n}^{(i)}(v_n)-E[U_{r_n-l_n}^{(i)}(v_n)]\right)+
a\sum_{i=1}^{k_n}\left(V_{l_n}^{(i)}(v_n)-E[V_{l_n}^{(i)}(v_n)]\right)\nonumber\\[0.3cm]
&&+ b\sum_{i=1}^{k_n}\left(\indi_{\{U_{r_n-l_n}^{(i)}(v_n)>0\}}-
E[\indi_{\{U_{r_n-l_n}^{(i)}(v_n)>0\}}]\right)+\nonumber\\[0.3cm]
&&+ b\sum_{i=1}^{k_n}\left(\indi_{\{U_{r_n-l_n}^{(i)}(v_n)>0,\ U_{r_n-l_n}^{(i)}(v_n)=0\}}-
E[\indi_{\{U_{r_n-l_n}^{(i)}(v_n)>0,\ U_{r_n-l_n}^{(i)}(v_n)=0\}}]\right).\label{eq3}
\end{eqnarray}
Now, since (\ref{eq3}) holds, to show (\ref{conv_normal2}) we have to prove the following
\begin{eqnarray}
&&c_n^{-1/2}
\left(a\displaystyle{\sum_{i=1}^{k_n}\left(U_{r_n-l_n}^{(i)}(v_n)-E[U_{r_n-l_n}^{(i)}(v_n)]
\right)+b
\sum_{i=1}^{k_n}\left(\indi_{\{U_{r_n-l_n}^{(i)}(v_n)>0\}}-
E[\indi_{\{U_{r_n-l_n}^{(i)}(v_n)>0\}}]\right)}\right)\nonumber\\
&&\hspace{2cm}\xrightarrow [n\rain]{d} \mathcal{N}
(0,\ \nu(a^2\eta \sigma^2+2ab+b^2\eta)),\label{conv_normal3}\\[0.4cm]
&&c_n^{-1/2} \sum_{i=1}^{k_n}(V_{l_n}^{(i)}(v_n)-E[V_{l_n}^{(i)}(v_n)])\xrightarrow [n\rain]{P}0,\label{eq5}\\[0.4cm]
&&c_n^{-1/2} \sum_{i=1}^{k_n}\left(\indi_{\{N_{r_n}^{(i)}>0,\ U_{r_n-l_n}^{(i)}(v_n)=0\}}-E[\indi_{\{N_{r_n}^{(i)}>0,\ U_{r_n-l_n}^{(i)}(v_n)=0\}}]\right)\xrightarrow [n\rain]{P}0.\hspace{1.5cm}\label{eq6}
\end{eqnarray}

Lets first prove (\ref{conv_normal3}). The summands in $\sum_{i=1}^{k_n}(aU_{r_n-l_n}^{(i)}(v_n)+b\indi_{\{U_{r_n-l_n}^{(i)}(v_n)>0\}})$ are functions of indicator variables that are at least $l_n$ time units apart from each other, therefore for each $t\in \R$
\begin{eqnarray}
\lefteqn{\left|E\left[\exp\left( {\rm{i}}tc_n^{-1/2}\sum_{i=1}^{k_n}\left(aU_{r_n-l_n}^{(i)}(v_n)+
b\indi_{\{U_{r_n-l_n}^{(i)}(v_n)>0\}}\right)\right)\right]-\right.}\\[0.3cm]
&&-\left.\prod_{i=1}^{k_n}E\left[{\rm{i}}tc_n^{-1/2}\left(aU_{r_n-l_n}^{(i)}(v_n)+
b\indi_{\{U_{r_n-l_n}^{(i)}(v_n)>0\}}\right)\right]\right|\leq 16k_n\alpha_{n,l_n},\label{eq7}
\end{eqnarray}
where, $\rm{i}$ is the imaginary unit, from repeatedly using a result in Volkonski and Rozanov (1959) \cite{ro} and the triangular inequality. Now, since condition $\Delta(v_n)$ holds for ${\bf{X}}$  (\ref{eq7}) tends to zero and so we can assume that the summands are i.i.d.. Therefore, in order to apply Lindberg's Central Limit Theorem we need to verify that
\begin{equation}
\frac{k_n}{c_n} E\left[\left(aU_{r_n-l_n}(v_n)+b\indi_{\{U_{r_n-l_n}(v_n)>0\}}\right)^2\right]\xrightarrow [n\rain]{}\nu(a^2\eta \sigma^2+2ab+b^2\eta),\label{eq8}
\end{equation}
since $\frac{k_n}{c_n}E^2[U_{r_n-l_n}(v_n)]\xrightarrow [n\rain]{}0,$ and Lindberg's Condition
\begin{equation}
\frac{k_n}{c_n}E\left[\left(aU_{r_n-l_n}(v_n)+b\indi_{\{U_{r_n-l_n}(v_n)>0\}}\right)^2
\indi_{\{(aU_{r_n-l_n}(v_n)+b\indi_{\{U_{r_n-l_n}(v_n)>0\}})^2>\epsilon c_n\}}\right]
\xrightarrow [n\rain]{}0,\quad \textrm{for all } \epsilon>0,\label{eq9}
\end{equation}
with $U_{r_n-l_n}(v_n)=U_{r_n-l_n}^{(1)}(v_n).$

From the definition of $\widetilde{N}_{r_n}^{(i)}(v_n)$ and $V_{l_n}^{(i)}(v_n)$ and assumption (2) we have that
\begin{equation}
\frac{k_n}{c_n}E[V_{r_n}^2(v_n)]\leq \left(\frac{r_n}{l_n}\right)^{-1}\frac{k_n}{c_n}E[\widetilde{N}_{r_n}^2(v_n)]\xrightarrow [n\rain]{}0,\label{asterisco}
\end{equation}
with $V_{r_n}(v_n)=V_{r_n}^{(1)}(v_n).$ Now, by Cauchy-Schwarz's inequality
$$\frac{k_n}{c_n}E[\widetilde{N}_{r_n}(v_n)V_{r_n}(v_n)]\leq \sqrt{E[\widetilde{N}_{r_n}^2(v_n)]E[V_{r_n}^2(v_n)]}\xrightarrow [n\rain]{}0,$$
thus
\begin{equation}
\frac{k_n}{c_n}E[U_{r_n-l_n}^2(v_n)]=\frac{k_n}{c_n}
E[(\widetilde{N}_{r_n}(v_n)-V_{r_n}(v_n))^2]\xrightarrow [n\rain]{}\eta \nu \sigma^2.\label{eq10}
\end{equation}
On the other hand, since $k_n(r_n-l_n)\sim n,$ Theorem 2.1 implies that
\begin{equation}
\frac{k_n}{c_n}E\left[\indi_{\{U_{r_n-l_n}(v_n)>0\}}\right]=\frac{k_n}{c_n}
P(U_{r_n-l_n}(v_n)>0)\xrightarrow [n\rain]{}\eta \nu.\label{eq11}
\end{equation}
Furthermore, by definition (\ref{nivel})
\begin{equation}
\frac{k_n}{c_n}E\left[U_{r_n-l_n}(v_n)\indi_{\{U_{r_n-l_n}(v_n)>0\}}\right]=
\frac{k_n}{c_n}E\left[U_{r_n-l_n}(v_n)\right]=\frac{k_n(r_n-l_n)}{c_n}P(X_1\leq v_n<X_2)\xrightarrow [n\rain]{}\nu.\label{eq12}
\end{equation}
(\ref{eq10})-(\ref{eq12}) prove (\ref{eq8}) and since Lindberg's Condition follows immediately from assumption (1), (\ref{conv_normal3}) is proven.

Finally, since Theorem 2.1 of Hsing (1991) \cite{hsing1} holds for $\widetilde{N}_{r_n}(v_n),$ it implies (\ref{eq5}) and (\ref{eq6}) because $\frac{k_n}{c_n}E[V_{l_n}^2(v_n)]\xrightarrow [n\rain]{}0$ by (\ref{asterisco}) and $$\frac{k_n}{c_n}E[\indi^2_{\{\widetilde{N}_{r_n}^{(i)}(v_n)>0,\ U_{r_n-l_n}^{(i)}(v_n)=0\}}]=\frac{k_n}{c_n}(P(\widetilde{N}_{r_n}^{(i)}(v_n)>0)-
P(U_{r_n-l_n}^{(i)}(v_n))>0)\xrightarrow [n\rain]{}0$$ by Theorem 2.1 and (\ref{eq11}). This concludes the proof.\hfill $\square$

\subsection{Proof of Corollary \ref{cor:normZ}}

Since $\eta_n\xrightarrow [n\rain]{}\eta,$ by Theorem 2.1, $\frac{k_n}{c_n}E[\widetilde{N}_{r_n}(v_n)]\xrightarrow [n\rain]{}\nu,$ and $c_n^{-1}\sum_{i=1}^{k_n}(\widetilde{N}_{r_n}^{(i)}(v_n)
-E[\widetilde{N}_{r_n}^{(i)}(v_n)])\xrightarrow [n\rain]{P}0,$ by Theorem 2.1 of Hsing (1991) \cite{hsing1} which holds for $\widetilde{N}{r_n}^{(i)}(v_n),$ the result now follows from the fact that
\begin{eqnarray*}
\lefteqn{\hspace{-1cm}\sqrt{c_n}(\widehat{\eta}_n^B-\eta_n)=\frac{1}{c_n^{-1}\sum_{i=1}^{k_n}
\widetilde{N}_{r_n}^{(i)}(v_n)}\left(c_n^{-1/2}\sum_{i=1}^{k_n}
\left(\indi_{\{\widetilde{N}_{r_n}^{(i)}(v_n)>0\}}-E[\indi_{\{\widetilde{N}_{r_n}^{(i)}(v_n)>0\}}]
\right)\right.-}\\[0.3cm]
&&\hspace{5.5cm}\left.-\eta_n c_n^{-1/2}\sum_{i=1}^{k_n}\left(\widetilde{N}_{r_n}^{(i)}(v_n)-
E[\widetilde{N}_{r_n}^{(i)}(v_n)]\right)\right)
\end{eqnarray*} and Theorem 2.3.  \hfill $\square$

\subsection{Proof of Theorem \ref{teo:aleatorio}}

Since (\ref{p1}) holds, we only need to show that
\begin{equation}
\sum_{j=1}^{+\infty}(\widetilde{\pi}_n(j;\widehat{v}_n)-\widetilde{\pi}_n(j;v_n))\xrightarrow [n\rain]{P}0. \label{pis}
\end{equation}
Lets start by noting that for $v_n^{(\tau)}$ such that $P(X_1>v_n^{(\tau)})\sim c_n\tau/n$ as $n\to +\infty$ and $\epsilon >0$ we have
\begin{eqnarray}
\lefteqn{\hspace{-1.2cm}\lim_{\epsilon\to 0}\lim_{n\to +\infty} \frac{k_n}{c_n}\left|E\left[\widetilde{N}_{r_n}(v_n^{(\tau+\epsilon)})\indi_{\{\widetilde{N}_{r_n}(v_n^{(\tau+\epsilon)})>0\}}
\right]-E\left[\widetilde{N}_{r_n}(v_n^{(\tau)})\indi_{\{\widetilde{N}_{r_n}(v_n^{(\tau)})>0\}}
\right]\right|}\nonumber\\[0.3cm]
&=&\lim_{\epsilon\to 0}\lim_{n\to +\infty}\frac{k_nr_n}{c_n}(P(X_1\leq v_n^{(\tau+\epsilon)}<X_2)-P(X_1\leq v_n^{(\tau)}<X_2))\nonumber\\[0.3cm]
&=&\lim_{\epsilon\to 0}\lim_{n\to +\infty}\frac{k_nr_n}{c_n}(P(X_1> v_n^{(\tau+\epsilon)})-P(X_1> v_n^{(\tau)}))=0.\label{pi_cov}
\end{eqnarray}
(\ref{pi_cov}) proves condition b) of Theorem 2.3 in Hsing (1991) \cite{hsing1} which holds for $\widetilde{N}_{r_n}(v_n),$ where $T(j)=j\indi_{\{j\geq 1\}}.$ The other conditions have been verified in the proof of Theorem 2.2 as well as the conditions of Corollary 2.4 in Hsing (1991) \cite{hsing1} for $\widetilde{N}_{r_n}(v_n).$ Therefore (\ref{pis}) holds, completing the proof.\hfill $\square$
\end{appendix}\vspace{0.3cm}

\noindent {\bf{Acknowledgements}} We acknowledge the support from research unit ``Centro de Matemática e Aplicações'' of the University of Beira Interior, through the research project UID/MAT/00212/2013.

\end{document}